\documentclass[11pt,reqno]{amsart}

\usepackage{marginnote,amssymb,mathrsfs,latexsym,verbatim,pdfsync,color}

\newcommand{\bbC}{{\mathbb C}}

\newcommand{\bbH}{{\mathbb H}} 

\newcommand{\bbR}{{\mathbb R}}

\def\bU{{\mathbf U}}

\def\cC{{\mathcal C}}
\def\cD{{\mathcal D}}

\def\cF{{\mathcal F}}

\def\cH{{\mathcal H}}

\def\cL{{\mathcal L}}

\def\cU{{\mathcal U}}

\def\sR{{\mathscr R}}

\def\Re{\operatorname{Re}}
\def\Im{\operatorname{Im}}

\def\la{\langle}
\def\ra{\rangle}

\def\eps{\varepsilon}
\def\z{\zeta} 

\def\ov{\overline}
\def\p{\partial}

\def\ms{\medskip}

\def\Hol{\operatorname{Hol}}
\def\tr{\operatorname{tr}}
\def\DA{\operatorname{DA}}
\def\HS{\operatorname{HS}}
\def\ran{\operatorname{ran}}

\def\i{\mathbf i}

\def\half{\textstyle{\frac12}}

\def\CB{\color{black} }

\newtheorem{thm}{Theorem}[section]
\newtheorem{prop}[thm]{Proposition}
\newtheorem{cor}[thm]{Corollary}
\newtheorem{lem}[thm]{Lemma}

\newtheorem{remark}[thm]{Remark}

\newtheorem{theorem}{Theorem}
\newtheorem*{Theorem}{Theorem}

\newtheorem*{definition}{Definition}

\begin{document}

\title[Paley--Wiener theorems 
on the Siegel upper half-space]{
Paley--Wiener theorems on the Siegel upper half-space} 
\author[N. Arcozzi, A. Monguzzi, M. M. Peloso, M. Salvatori]{Nicola
  Arcozzi, 
Alessandro Monguzzi, Marco M. Peloso, Maura Salvatori}
\address{Dipartimento di Matematica, Alma Mater Studorium Universit\`a di
  Bologna, Piazza di Porta San Donato 5, 40126 Bologna, Italy}
\address{Dipartimento di Matematica, Universit\`a degli Studi di
  Milano, Via C. Saldini 50, 20133 Milano, Italy}
\email{{\tt nicola.arcozzi@unibo.it}}
\email{{\tt alessandro.monguzzi@unimi.it}}
\email{{\tt marco.peloso@unimi.it}}
\email{{\tt maura.salvatori@unimi.it}}
\keywords{Siegel upper half-space, holomorphic function spaces, reproducing 
kernel Hilbert space, Drury--Arveson, Dirichlet, Hardy, Bergman spaces.}
\thanks{{\em Math Subject Classification} 30H99, 46E22, 30C15, 30C40.}
\thanks{Authors are partially supported by the 2015 PRIN grant
  \emph{Real and Complex Manifolds: Geometry, Topology and Harmonic Analysis}  
  of the Italian Ministry of Education (MIUR)}

\begin{abstract}
In this paper we study spaces of 
holomorphic  functions on the Siegel upper half-space $\cU$ and
prove Paley--Wiener type theorems for such spaces.
The 
boundary of $\cU$ 
can be identified with the Heisenberg group $\bbH_n$.   Using the
group
Fourier transform on $\bbH_n$, 
Ogden--Vagi \cite{OV} proved a Paley--Wiener theorem for the
Hardy space $H^2(\cU)$.

We consider a scale of Hilbert spaces on  
$\cU$ that includes the Hardy space, the weighted Bergman spaces,
the weighted Dirichlet spaces, and in particular the Drury--Arveson
space, 
and the Dirichlet space $\cD$.  
For each of these spaces, we 
prove a Paley--Wiener theorem, some structure
theorems,
and provide some applications.  In particular we prove that the norm of the
Dirichlet space modulo constants $\dot{\cD}$ is the unique Hilbert
space norm that is invariant under the action of the group of automorphisms of
$\cU$. 
\end{abstract}
\maketitle

\section{Introduction and statement of the main results}

Let $\bbC_+$ be the upper half-plane $\{z=x+iy\in\bbC:\, y>0\}$.
Let $H^2(\bbC_+)$ denote the Hardy space, that is the space of holomorphic
functions in $\bbC_+$ such that
$$
\| f\|_{H^2(\bbC_+)}^2 :=
\sup_{y>0} \int_{-\infty}^{+\infty} |f(x+iy)|^2\, dx <+\infty \,.
$$
The classical Paley--Wiener theorem \cite{PW} says that, given $f\in
H^2(\bbC_+)$ there exists $g\in L^2(0,+\infty)$ such that
\begin{equation}\label{PW-eq-n=0;1}
f(z) = \frac{1}{2\pi} \int_0^{+\infty} e^{iz\xi} g(\xi)\, d\xi
\,,
\end{equation}
and 
\begin{equation}\label{PW-eq-n=0;2}
\|f\|_{H^2(\bbC_+)}^2 = \frac{1}{2\pi} \|g\|_{L^2(0,+\infty)}^2 \,.
\end{equation}
Conversely, given any $g\in L^2(0,+\infty)$, defining $f$ as in
\eqref{PW-eq-n=0;1}, we have that $f\in
H^2(\bbC_+)$ and \eqref{PW-eq-n=0;2} holds.     Since then, the 
Fourier transform on the real line has appeared as a fundamental tool
in modern complex analysis and in the theory of holomorphic function
spaces in $\bbC_+$.
  We mention, for instance,
the regularity 
of projection operators,
 the  boundary behavior and  growth conditions
 of holomorphic functions,
the  boundedness and compactness of Hankel and
 Toeplitz operators, just to name some of the most important, see
 e.g. \cite{Poltoratski,Seip}.
The Paley--Wiener theorem has been extended to other Hilbert
function spaces on $\bbC_+$: the weighted Bergman spaces (see
e.g. \cite{FK,camerun,DG}), and more recently to the Dirichlet space
(\cite{Kucik,Kucik2}). 

The
 domain $\bbC_+$ is biholomorphic equivalent to the unit disk in  the
 plane.   In principle, it is possible to transfer analogous results
 from the unit disk to $\bbC_+$.  However, often, it is far more natural to
 study a problem directly on the unbounded domain $\bbC_+$.

In this paper we wish to extend the approach described above to
$\bbC^{n+1}$, where we always assume that $n\ge1$.  

The  Siegel upper-half space  is the domain in
$\bbC^{n+1}$ 
$$
\cU =\bigl\{ \z=(\z',\z_{n+1})\in\bbC^n\times\bbC:\, \Im
\z_{n+1} >\textstyle{\frac14} |\z'|^2 \bigr\} \,,
$$
and we denote by $\varrho(\z)=
\Im \z_{n+1}-\frac14|\z'|^2$ its defining function.  
The domain $\cU$ is biholomorphic to
the unit ball $B$ in $\bbC^{n+1}$ via the (multi-dimensional) Cayley
transform $\cC:B\to\cU$,
$$
\cC(\omega)=\Bigl( \frac{2\omega'}{1-\omega_{n+1}}, 
i\frac{1+\omega_{n+1}}{1-\omega_{n+1}}\Bigr)\,. 
$$ 
  The boundary $\p\cU$ of $\cU$ can be endowed with the structure of a
nilpotent Lie group, namely the Heisenberg group $\bbH_n$.   Thus, it
is possible to use the group Fourier transform on $\p\cU$ to
characterize the boundary values of holomorphic functions on $\cU$
in various Hilbert function spaces, and therefore to take a first step in the
program outlined in the case of $\bbC_+$. In particular we prove
Paley--Wiener type theorems for functions in weighted Bergman and
 Dirichlet
spaces, and on the Dirichlet space. In fact, we show that the latter one is the unique Hilbert space modulo constants that is invariant under the group of automorphisms of $\cU$.
We now describe  the
content of the present paper in greater
details. \ms

The boundary $\p\cU$ is characterized by the points of
$\bbC^{n+1}$ such that $\varrho(\z)=0$, that is,
$\z=(\z',t+\textstyle\frac{i}{4}|\z'|^2)$, with $t\in\bbR$.  
We
introduce 
a parametrization of $\cU$ by means of a foliation of copies
of the boundary.  
We set $\bU= \bbC^n\times\bbR\times(0,+\infty)$.
Given $\z=(\z',\z_{n+1})\in \cU$, we define $\Psi(\z',\z_{n+1})=(z,t,h) \in
\bU$ by
\begin{equation}\label{h-par}
\begin{cases}
z=\z'\cr
t= \Re \z_{n+1}\cr
h =\Im \z_{n+1} -\frac14|\z'|^2 \,.
\end{cases}
\end{equation}
Then,  
$\Psi: \ov\cU\to
\ov\bU$ is a $C^\infty$-diffeomorphism,  
and $\Psi^{-1}$
is given by
\begin{equation}\label{h-coordinates}
\Psi^{-1} (z,t,h) = \bigl(z,t+i\textstyle{\frac14}|z|^2+ih\bigr) =:  (\z',\z_{n+1})
\,. 
\end{equation}
Notice that
$h=\varrho(\z',\z_{n+1})$.   
When $h=0$,  we write $[z,t]$ in
place of $(z,t,0)$.  
The points on the boundary act on $\overline{\cU}$ as biholomorphic
maps in the following way.  
For $[z,t]\in\p\cU$, we define 
\begin{equation}\label{Phi}
\Phi_{[z,t]} (\omega',\omega_{n+1})
= \bigl(\omega'+z, \omega_{n+1} +t
+i\textstyle{\frac14}|z|^2+ \frac i2 \omega'\cdot\bar z \bigr) \,,
\end{equation}
where $\omega'\cdot \bar z=\sum_{j=1}^n\omega_j\bar z_j$ denotes the
hermitian inner product in $\bbC^n$. 
Notice that
\begin{align*}
\varrho\Bigl(
\Phi_{[z,t]}
(\omega',\omega_{n+1})\Bigr) 
& = \varrho(\omega',\omega_{n+1})\,,  
\end{align*}
that is, the maps $\Phi_{[z,t]}$ preserve the
defining function $\varrho$. 
 In particular, for $(\omega',\omega_{n+1})\in\p\cU$ and
$[w,s]=\Psi(\omega',\omega_{n+1})$, by \eqref{Phi} we have
\begin{align}
\Phi_{[z,t]} \big(
(\omega',\omega_{n+1})\big) & = \Phi_{[z,t]} \Big( \Psi^{-1}
(w,s,0)\Big) = 
\Phi_{[z,t]} \bigl(w,s+i\textstyle{\frac14}|w|^2\bigr) \notag \\
 & =\bigl(w+z,s+\textstyle{\frac i4}|w|^2+t+ \textstyle{\frac i4}|z|^2
+\textstyle{\frac i2} w\cdot\bar
z\bigr) \notag \\
& =\bigl[ w+z,s+t -\half\Im(w \cdot\bar
z)\bigr] \notag \\ 
& =: [w,s][z,t] \, .  
\label{product-on-bUn}
\end{align}
Therefore, it is possible to introduce a group structure on $\p\cU$ itself.
\begin{definition}{\rm
The Heisenberg group $\bbH_n$  is the set $\bbC^n\times\bbR$ endowed
with product 
$$
[w,s][z,t]  =  \big[ w+z, s+t -\textstyle{\frac12} \Im (w\cdot\bar z)\big] \,.
$$
}
\end{definition}

The Heisenberg group $\bbH_n$ is a nilpotent Lie group of step 2, and
the Lebesgue measure on $\bbC^n\times\bbR$ coincides with both the
right and left Haar
measure on $\bbH_n$.  In other words, the Lebesgue measure is both
right and left translation invariant. 

If $x$ is a vector of the Euclidean space $\bbR^d$, we denote by $dx$
the Lebesgue measure in $\bbR^d$.
Notice that, since $|\text{det Jac}\Psi|=1$, for $F$ integrable on
$\cU$, setting $\widetilde F=F\circ\Psi^{-1}$ and $\widetilde
F_h[z,t]:=\widetilde F(z,t,h)$,
we have
$$
\int_\cU F(\z)\, d\z = 
\int_\bU \widetilde F (z,t,h)\,  dzdtdh =
\int_0^{+\infty} \int_{\bbH_n}
\widetilde F_h [z,t]\, dzdt dh \,.
$$

We now introduce the Hilbert function
spaces object of our study.

\begin{definition}\label{Hilber-fnc-sp-def}{\rm  
For $\nu>-1$, we consider the {\em  weighted Bergman spaces}
$A^2_\nu $
\begin{equation*}
A^2_\nu  = \Big\{ F\in\Hol(\cU): \,  \|F\|_{A^2_\nu}^2 :=\int_\cU |F(\z)|^2 \rho(\z)^\nu\,
d\z = 
\int_\bU 
|\widetilde F (z,t,h)|^2 \, h^\nu\, dzdt dh
<+\infty\Big\} \,.
\end{equation*}

For $-n-2<\nu<-1$ and $m$ a positive integer such that
$2m+\nu>-1$, the {\em weighted Dirichlet spaces}
 are defined as follows
\begin{equation}\label{w-Dir-space-def}
\begin{aligned} 
& \cD_{\nu,(m)}
=\Big\{ 
F\in\Hol(\cU):&
{\rm (i)}\ & \displaystyle { \lim_{|\z'|\le R,\, \Im\z_{n+1}\to+\infty} 
  F(\z) =0} \,;& \smallskip\\
& &  {\rm (ii)} &  \displaystyle {\int_\cU |\rho^m(\z) \p_{\z_{n+1}}^m F(\z)|^2 \,
\rho^\nu(\z) d\z <+\infty  }\, \Big\}\,.
\end{aligned}
\end{equation}
For $F$ as above, we define the norm on $\cD_{\nu,(m)}$ as
$$
\|F\|_{\cD_{\nu,(m)}}^2 =
\int_\cU |\rho^m(\z) \p_{\z_{n+1}}^m F(\z)|^2 \,
\rho^\nu(\z) d\z  \,.
$$

Finally, for $\nu=-n-2$ and $2m>n+1$, we define the {\em Dirichlet space} $\cD_{(m)}$ as
\begin{equation}\label{Dir-space-def}
\begin{aligned} 
& \cD_{(m)}
=\Big\{ 
F\in\Hol(\cU):&
{\rm (i)}\ 
& \displaystyle { \lim_{|\z'|\le R,\, \Im\z_{n+1}\to+\infty} \p_{\z_j}
  F(\z) =0}\ \text{for}\ j=1,\ldots,n+1  ;& \smallskip\\ 
& &  {\rm (ii)} &  \displaystyle {\int_\cU |\rho^m(\z) \p_{\z_{n+1}}^m F(\z)|^2 \,
\rho^{-n-2}(\z) d\z <+\infty  }\, \Big\}\,,
\end{aligned}
\end{equation}
with norm given by
\begin{equation}\label{Dir-norm-def}
\| F \|_{\cD_{(m)}}^2 = \|  \p_{\z_{n+1}}^m F
\|_{A^2_{2m-n-2}}^2 + |F(\i)|^2
\end{equation}
where $\i=(0',i)\in\cU$.
}\end{definition}
\ms

The case $\nu=-1$ corresponds to the classical  {\em Hardy space} $H^2$, defined as 
\begin{equation}\label{Hardy-space}
H^2  = \Big\{ F\in\Hol(\cU): \,  
\|F\|_{H^2}^2 :=\sup_{h>0}\int_{\bbH_n} |\widetilde F_h[z,t]|^2 \,
dzdt <+\infty\Big\} \,.
\end{equation}

We point out that when $2m+\nu>-1$, the spaces $\cD_{\nu,(m)}$  all
coincide, with the same norms up to a positive constant multiple (see 
Theorem \ref{PW-thm-Dir}).  Thus, when the choice of the norm
is unambiguous,
we simply denote them by
$\cD_\nu$.
Analogously, 
 the spaces $\cD_{(m)}$ do not depend on the choice of the
 integer $m$ (see Theorem \ref{main-4}), and we denote them
by $\cD$.  
  Moreover, we will show
that the norm of  $\cD$ modulo constants is invariant under the
automorphism group; hence 
it is legitimate to call the space $\cD$ the {\em Dirichlet space} on the
Siegel upper half-space. 

\ms
Our main technical tool is the Fourier transform on the Heisenberg
group.  For this, and other basic facts concerning the Heisenberg
group, we refer the reader to \cite{Folland} and \cite{R}.

Let $\lambda\in\bbR\setminus\{0\}$.   We 
set
\begin{equation}\label{Fock-sp}
\cF^\lambda 
=\bigg\{ F \in\Hol(\bbC^n):\ 
\bigg( \frac{|\lambda|}{2\pi} \bigg)^n \int_{\bbC^n} |F(z)|^2\,
e^{-\frac\lambda2 |z|^2}dz <+\infty \bigg\}
\end{equation}
when $\lambda>0$,   and $\cF^\lambda = \cF^{|\lambda|}$ when  
$\lambda<0$, and call this space the  Fock space.  We present
further properties of such space in Subsection \ref{Fock-space-subsec}.

For $\lambda\in\bbR\setminus\{0\}$ and $[z,t]\in\bbH_n$, the
Bargmann 
representation $\sigma_\lambda[z,t]$ is the operator acting on
$\cF^\lambda$ given by, 
\begin{equation}\label{Barg-repr-lambda-pos} 
\sigma_\lambda[z,t] F (w)
= e^{i\lambda t-\frac\lambda2 w\cdot\ov z -\frac\lambda4 |z|^2}
F(w+z) 
\end{equation}
if $\lambda>0$, and, if $\lambda<0$, as 
$\sigma_\lambda[z,t] = \sigma_{-\lambda} [\ov z,-t]$, that is, 
\begin{equation}\label{Barg-repr-lambda-neg} 
\sigma_\lambda[z,t] F (w)
= e^{i\lambda t+\frac\lambda2 w\cdot z +\frac\lambda4 |z|^2}
F(w+\bar z) \,.
\end{equation}

If $f\in L^1(\bbH_n)$, for $\lambda\in\bbR\setminus\{0\}$, 
$\sigma_\lambda(f)$ is the operator 
acting on $\cF^\lambda$ as
$$
 \sigma_\lambda (f) F (w) 
= \int_{\bbH_n} f[z,t] \sigma_\lambda[z,t]F (w)\, dzdt\,. 
$$

Before stating our main results, we recall a result proved by
Ogden and Vagi \cite{OV}, that extends the classical Paley--Wiener theorem for
the Hardy space,
from the upper half-plane to the case of 
$\cU$.  We point out however, that Ogden and Vagi proved their main
result in the case of   Siegel domains of type $I\!I$.  It
would certainly be  of interest to extend our results to the latter
more general class of domains. 

\ms

\begin{Theorem}[\cite{OV}]
Let $F\in H^2$. 
Then, there exists $\widetilde F_0\in L^2(\bbH_n)$
such that $ \widetilde F_h\to\widetilde F_0$  in $L^2(\bbH_n)$, as $h\to 0^+$. 
Moreover,
the function $\widetilde F_0$ is such that
\begin{itemize}
\item[(i)] $\|F\|_{H^2}  = \|\widetilde F_0\|_{L^2(\bbH_n)}$;\smallskip
\item[(ii)]  $\sigma_\lambda(\widetilde F_0) =0$ when
  $\lambda>0$;\smallskip
\item[(iii)] $\ran \big(\sigma_\lambda(\widetilde F_0)\big) \subseteq
  \operatorname{span}\{ 1\}$ for $\lambda<0$.
\end{itemize}

Conversely, if $f\in L^2(\bbH_n)$ is such that (ii) and (iii) are
satisfied, then setting
\begin{equation}\label{ext-form}
F(\z) = \widetilde F_h[z,t] =
\frac{1}{(2\pi)^{n+1}} \int_{-\infty}^0 e^{h\lambda} \tr \big(
\sigma_\lambda(f) \sigma_\lambda[z,t]^*\big)\, |\lambda|^n d\lambda\,,
\end{equation}
then
 $F\in H^2$ is such that $\widetilde F_0=f$ and {\rm (i)-(iii)} hold.
\end{Theorem}

We also need the following
\begin{definition}\label{F-T-A2nu}{\rm
For $\nu\in\bbR$ we define  the space 
$\cL^2_\nu $ as the space of functions $\tau$
 on $\bbR\setminus\{0\}$ such that:
\begin{align}\label{L2-nu}
 \begin{split}
\rm{(i)}&\textrm{ $\tau(\lambda) \in \HS(\cF^\lambda)$  for
  every $\lambda$, i.e., $\tau(\lambda):\cF^\lambda\to\cF^\lambda$ is a Hilbert--Schmidt operator;}\\
\rm{(ii)}&\textrm{ $\tau(\lambda)=0$ for
  $\lambda>0$;}\\
 \rm{(iii)}& \textrm{ $\ran( \tau(\lambda)) \subseteq
   \operatorname{span}\{1\}$;}\\
\rm{(iv)}& \textrm{ $
\displaystyle{ \|\tau\|_{\cL^2_\nu}^2: = \frac{1}{(2\pi)^{n+1}} 
\int_{-\infty}^0 \|\tau(\lambda)\|_{\HS}^2 \,
|\lambda|^{n-(\nu+1)} d\lambda <+\infty }$,  where 
$\|\cdot\|_{\HS} :=\|\cdot\|_{\HS(\cF^\lambda)}$.} 
 \end{split}
\end{align}}
\end{definition}

Our first main result is the following Paley--Wiener type theorem for
the weighted Bergman spaces $A^2_\nu$.

\begin{theorem}\label{PW-thm-Berg} 
Let $\nu>-1$ be fixed. Given $F\in A^2_\nu$, there exists $\tau\in \cL^2_\nu $ such
that, for $\z\in\cU$, 
\begin{equation}\label{PW-A2nu-eq1}
F(\z) = \widetilde F_h[z,t] =
\frac{1}{(2\pi)^{n+1}} \int_{-\infty}^0 
 e^{h\lambda} \tr \big(
\tau(\lambda)\sigma_\lambda[z,t]^*\big) \, |\lambda|^n d\lambda\,,
\end{equation}
and
\begin{equation}\label{PW-A2nu-eq2}
\| F\|_{A^2_\nu}^2
= \frac{\Gamma(\nu+1)}{2^{\nu+1}} \| \tau\|_{\cL^2_\nu}^2 \,. 
\end{equation}

Conversely, given $\tau\in\cL^2_\nu$, let $F$ be defined as in
\eqref{PW-A2nu-eq1}. Then $F\in A^2_\nu$ and \eqref{PW-A2nu-eq2}
holds. 
\end{theorem}

Next we consider the case of weighted Dirichlet spaces.

\begin{theorem}\label{PW-thm-Dir}
Let $-(n+2)< \nu<-1$, and let $m>-\frac{\nu+1}{2}$. 
Let $F\in \cD_{\nu,(m)}$.  Then, there exists $\tau\in \cL^2_\nu $ such
that, for $\z\in\cU$,
\begin{equation}\label{PW-Dnu-eq1}
F(\z) = \widetilde F_h[z,t] =
\frac{1}{(2\pi)^{n+1}} \int_{-\infty}^0 
 e^{h\lambda} \tr \big(
\tau(\lambda)\sigma_\lambda[z,t]^*\big) \, |\lambda|^n d\lambda\,,
\end{equation}
and
\begin{equation}\label{PW-Dnu-eq2}
\| F\|_{\cD_{\nu,(m)}}^2
= \frac{\Gamma(2m+\nu+1)}{2^{2m+\nu+1}} 
\| \tau\|_{\cL^2_\nu}^2 \,. 
\end{equation}

Conversely, given $\tau\in\cL^2_\nu$, let $F$ be defined as in
\eqref{PW-Dnu-eq1}. Then $F\in \cD_{\nu,(m)}$ and \eqref{PW-Dnu-eq2}
holds. 

Therefore, for each $m>-\frac{\nu+1}{2}$, the spaces $\cD_{\nu,(m)}$
all coincide and their norms satisfy \eqref{PW-Dnu-eq2}.  
\end{theorem}

Hence, if no confusion arises, 
we simply write $\cD_\nu$ in place of $\cD_{\nu,(m)}$.  \ms

In the case $\nu=-n-1$, $\cD_\nu$ is called the {\em Drury--Arveson
  space} and we denote it by $\DA$.   The 
Drury--Arverson
  space on the unit ball $B$ has drawn a great deal of interest
  in the recent years,
  see \cite{Arveson,Drury,ARS,T,ARSW-2variations,CSW,VW,RS}, and references
  therein, to name a few.  When
  $n\ge1$, $\DA$ plays a role similar to the one played by the Hardy space on the unit
  disk, 
and for this reason it is sometimes denoted as $H^2_{n+1}$.   If $f$ is
  holomorphic on $B$, $f(\z)=\sum_{|\alpha|\ge0}  a_\alpha
  \z^\alpha$, the norm in $\DA(B)$ is given by
$$
\| f\|_{\DA(B)}^2 
= \sum_{|\alpha|\ge0} \frac{\alpha!}{|\alpha|!} |a_\alpha|^2\,.
$$
However, to the best of our knowledge, no integral representation of
this norm has been found.   
In this paper we provide such a description, see Theorem \ref{DA-ball-thm}.\footnote{We thank H. Turgay Kaptanoğlu for informing us, on date May 30, 2024, that this formula and related ones already appeared in ``D. Alpay and T. Kaptanoğlu, \textit{Integral formulas for a sub-Hardy Hilbert space on the ball with complete Nevanlinna-Pick reproducing kernel}, C. R. Acad. Sci. Paris Sér. I Math.333(2001), no.4, 285–290.''.}

The last main result is the following.
\begin{theorem}\label{main-4} Let $m>\frac{n+1}{2}$ be fixed. Let
  $F\in \cD_{(m)}$.  Then, there exists $\tau\in\cL^2_{-n-2}$ such 
that, for $\z\in\cU$,
\begin{equation}\label{PW-D-eq1}
F(\z)= \widetilde F_h[z,t] =
\frac{1}{(2\pi)^{n+1}} \int_{-\infty}^0
\tr\Big(\tau(\lambda)\big(e^{\lambda
  h}\sigma_\lambda[z,t]^*-e^\lambda\sigma_\lambda(0,0)^*\big) \Big)\,
|\lambda|^n d\lambda + c\,, 
\end{equation}
where $c=F(\i)$, and
\begin{equation}\label{PW-D-eq2}
\| F\|_{\cD_{(m)}}^2
= \frac{\Gamma(2m-n-1)}{2^{2m-n-1}} 
\| \tau\|_{\cL^2_{-n-2}}^2+|F(\i)|^2 \,. 
\end{equation}

Conversely, given $\tau\in\cL^2_{-n-2}$, let $F$ be defined as in
\eqref{PW-D-eq1}. Then $F\in \cD_{(m)}$, $c= F(\i)$ and \eqref{PW-D-eq2}
holds.

Therefore, for each $m>\frac{n+1}{2}$, the spaces $\cD_{(m)}$
all coincide and their norms satisfy \eqref{PW-D-eq2}. 
\end{theorem}

Hence, if no confusion arises, 
we simply write $\cD$ in place of $\cD_{(m)}$.   We shall also denote
by $\dot{\cD}$ the quotient space $\cD/\bbC$, endowed with any of the
norms $\| \p_{\z_{n+1}}^m F\|_{A^2_{2m-n-2}}$.  We are going to show
  that $\dot{\cD}$ is the unique Hilbert space of functions modulo
  constants that is invariant under
composition with automorphisms, see Theorem \ref{Mobius-invariance}. We would like to point out that, when $\nu=-n-2$, even
given the integrability condition of the derivative of sufficiently
high order $m$, 
 it is not possible to find an  anti-derivative of order $m$ that vanishes as
 $\Im\z_{n+1}\to+\infty$. Hence, the decay property in
 \eqref{Dir-space-def} is required on the gradient of the function and
 not on the function itself. 
We will comment and make more remarks in 
Sections \ref{w-Dir-sec} and \ref{Dir-sec}. 
\ms

Beside their intrinsic interest, 
there are several reasons to  
study Paley--Wiener type theorems. 
All the spaces we are considering are Hilbert spaces, in particular
with 
a reproducing
kernel, and it is possible to define the same scale of space 
 with $p\neq2$. These spaces  are 
classical Besov--Sobolev spaces; for the case of
 the unit ball, see, e.g., \cite{BB,VW,Zhu}.  The boundary behavior
 of functions 
in the weighted Bergman spaces $A^p_\nu(\cU)$ with $\nu>-1$ was studied
by M. Feldman \cite{Feldman}, following the case of the upper half-plane in
$\bbC_+$ obtained by F. Ricci and M. Taibleson \cite{RT}.  
Among other results, we provide the 
 explicit expression of the reproducing kernels for all these spaces. 
These kernel are also the integral kernels for the corresponding
orthogonal projections.  It would be of interest to study the regularity
properties of such projections on the scale of the appropriate 
homogeneous Sobolev
spaces.
As potential application of our results, 
we also mention  the theory of invariant subspaces, in the
spirit of \cite{Lax}, \cite{DG}, e.g., that deal with this question in
the 1-dimensional setting of the upper half-plane. Furthermore, we
point out that the the invariance of the norm of the Dirichlet space 
under the composition with the automorphisms is much easier to
prove in the setting of the Siegel half-plane than in the unit
ball -- cfr. (2) in Theorem \ref{Mobius-invariance} and \cite[Theorem
6.15]{Zhu}. 
\CB  \ms

The paper is organized as follows. Section \ref{Sec-2} is a
preliminary section where we recall some standard results on the
Siegel half-space, the Heisenberg group and the Hardy space on
$\cU$. In Section \ref{Berg-section}, \ref{w-Dir-sec} and
\ref{Dir-sec} the weighted Bergman spaces, the weighted Dirichlet
spaces and the Dirichlet space are studied respectively. In Section
 \ref{remarks-sec} we provide the integral norm of the
Drury--Arveson 
space on the unit ball, and then we conclude  
with some final remarks and possible future
directions of research.

\section{Preliminaries}
\label{Sec-2} 

In this part we recall some well-known facts that will be used in what
follows.

\subsection{More on the Heisenberg group and the the Siegel upper
 half-space}

The following lemma is well known, see e.g. \cite[7.5.18]{Wallach}; we thank F. Ricci for pointing this reference to us.
 \begin{lem}\label{Aut-U}
The group $\operatorname{Aut}(\cU)$
of biholomorphic self-maps of $\cU$
is given by
$$
\operatorname{Aut}(\cU) = \bigcup_{\gamma\in\{ {\rm Id}, v\}} (MAN)\gamma(MAN)\,,
$$
where
\begin{itemize}
\item[(i)] $N=\big\{ \Phi_{[z,t]}:\, [z,t]\in\bbH_n\big\}$ (the
  subgroup of Heisenberg translations);\smallskip
\item[(ii)]$A=\big\{ D_\delta:\delta>0,\ 
D_\delta(\z',\z_{n+1})= (\delta\z',\delta^2\z_{n+1})\big\}$ (the
subgroup of non-isotropic dilations);\smallskip
\item[(iii)] $M=\big\{ U\in  U(n):\, 
\Phi_U(\z',\z_{n+1})= (U\z',\z_{n+1})\big\}$ (the subgroup of unitary
transformations in $\bbC^n$);\smallskip
\item[(iv)] $v(\zeta)=
\Bigl( \frac{i\z'}{\z_{n+1}}, -\frac{1}{\z_{n+1}}\Bigr)$ (the
inversion map).
\end{itemize}
\end{lem}

On $\bbH_n$
we define a {\em homogeneous}
norm 
by setting
\begin{equation*}
|[z,t]|_{\bbH_n} := 
\bigl(\textstyle{\frac{1}{16} } |z|^4 +t^2\bigr)^{1/4}\, .
\end{equation*}
This norm satisfies the following properties:
\begin{itemize}
\item[{\tiny$\bullet$}] $|[z,t]|_{\bbH_n}\ge0$ and it is 0 if and only if $[z,t]=(0,0)$;\smallskip
\item[{\tiny$\bullet$}] $\big|[z,t][w,s]\big|\le |[z,t]|+|[w,s]|$;\smallskip
\item[{\tiny$\bullet$}] $|D_\delta [z,t]|_{\bbH_n}=\delta |[z,t]|_{\bbH_n}$.
\end{itemize}
The topology induced by the metric $d_{\bbH_n} \bigl( [z,t],[w,s]\bigr)=
|[z,t][w,s]^{-1}|_{\bbH_n}$ is equivalent to the Euclidean topology of
$\bbC^n\times\bbR$. 
We also set
\begin{equation*}
B\big([z,t],r\big)
= \Big\{ [w,s]\in\bbH_n:\, |[w,s][z,t]^{-1}|_{\bbH_n} <r
\Big\}\,.
\end{equation*}
\ms

We recall that a holomorphic function $F$ 
satisfies the mean value property
$F(\z) = \frac{1}{|Q|} \int_Q F(\omega)\, d\omega$, where $Q=Q(\z,R)$ denotes
the polydisk $\{\omega:\, |\omega_j-\z_j|<R_j\}$ of polyradius $R$, 
contained in the region of 
holomorphy of $F$, and $|Q|$ is its Lebesgue measure. 
We will also consider a metric on $\cU$, which is somehow
conformally invariant.  If  $\Psi(\z',\z_{n+1})=(z,t,h)$ 
as in \eqref{h-par}, we set
$$
P ((z,t,h),r) 
= B\big([z,t],r\big)\times
\big\{ k:\,   |h-k|<r^2\big\} \quad \textrm{and} \quad P(\z,r)=\Psi^{-1}\big(P((z,t,h),r)\big).
$$
Then, we have 
\begin{itemize}
\item[{\tiny$\bullet$}] $|P (\z,r)| = c_n r^{2n+4}$;\smallskip
\item[{\tiny$\bullet$}] $\Phi_{[w,s]} \big(P (\z,r) \big) =
  P\big( \Phi_{[w,s]} (\z) ,r \big)$ ;\smallskip
\item[{\tiny$\bullet$}] $D_\delta \big(P (\z,r) \big) = P
  \big( D_\delta(\z) , r\delta\big)$.
\end{itemize}
It is elementary to see that a holomorphic function $F$ 
satisfies the submean 
value property
\begin{equation*}
|F(\z) | 
\le \frac{C}{|P(\z,r)|}\int_{|h-k|<r^2} \int_{B([z,t],r)} |\widetilde F(w,s,k)|\, dwds dk\,.
\end{equation*}
\ms

\subsection{The Fock space and the Fourier transform on the Heisenberg
  group}\label{Fock-space-subsec} 

Recall that the Fock space $\cF^\lambda$ is defined in
\eqref{Fock-sp} and thus has inner product 
$$
\la f,g\ra_{\cF^\lambda}
= \bigg( \frac{|\lambda|}{2\pi} \bigg)^n \int_{\bbC^n} f(z) \ov{g(z)}\,
e^{-\frac\lambda2 |z|^2}dz\,.
$$
Observe that 
$
\big(
\frac{|\lambda|}{2\pi} \big)^n e^{-\frac{|\lambda|}{2} |z|^2} dz$ 
is a probability
measure, and that the normalized monomials $\big\{
z^\alpha/\|z^\alpha\|_{\cF^\lambda}\big\}$, 
form a complete orthonormal
basis, and 
$$
\|z^\alpha\|_{\cF^\lambda}^2 
= \alpha!
\bigg(\frac{2}{|\lambda|}\bigg)^{|\alpha|}
\,.
$$
Moreover, $\cF^\lambda$ is a reproducing kernel Hilbert space,
with reproducing kernel $  e^{\frac{|\lambda|}{2} z\bar w} $,
\cite{Folland}. \ms

Introducing real coordinates on $\bbH_n$, $z_j=x_j+iy_j$,
$j=1,\dots,n$, then
$\bbH_n =\bbR^n\times \bbR^n\times\bbR$, 
and a basis for the left-invariant vector fields is $\{
X_1,\dots,X_n,Y_1,\dots,Y_n,T\}$, where
$$
X_j = \p_{x_j}-\textstyle{\frac{1}{2}}y_j \p_t\,,\ 
Y_j = \p_{y_j}+\textstyle{\frac{1}{2}}x_j \p_t\,,\ 
T= \p_t\,.
$$
A
basis for the complexified vector fields is $\big\{ Z_1,\dots,Z_n,\bar
Z_1,\dots,\bar Z_n,T\big\}$, where
$$
Z_j = \half(X_j-iY_j)= \p_{z_j} -\textstyle{\frac i4} \bar z_j \p_t,
\ 
\bar Z_j = \half(X_j+iY_j)= \p_{z_j} +\textstyle{\frac i4}  z_j \p_t, 
\  j=1,\dots,n\,, \ 
T\,,
$$
with non-trivial commutation rules 
$$
[Z_j,\bar Z_j]=\textstyle{\frac i2}T,\  j=1,\dots,n\,.
$$
We denote by $Z_j^{(R)}$ and $\bar Z_j^{(R)}$, resp., $j=1,\dots,n$,
the right-invariant 
vector fields that coincide with $Z_j $ and
$\bar Z_j$, resp., at the origin.  It turns out that
$$
Z_j ^{(R)}
= \p_{z_j} +\textstyle{\frac i4} \bar z_j \p_t,
\ 
\bar Z_j^{(R)} = \p_{\ov z_j}     -\textstyle{\frac i4} z_j \p_t
\,,  j=1,\dots,n\,.
$$
Then, the differentials of the Bargmann representations, that are
defined in \eqref{Barg-repr-lambda-pos}  and
\eqref{Barg-repr-lambda-neg}, can be computed 
to give, in particular:
\begin{itemize}
\item[(i)] for all $\lambda\neq0$, $d\sigma_\lambda(T) = i\lambda$;\smallskip
\item[(ii)] for $\lambda>0$,  $d\sigma_\lambda(\ov Z_j^{(R)}) = -\frac\lambda2 w_j$;\smallskip  
\item[(iii)]  for $\lambda<0$,   and
  $d\sigma_\lambda(\ov Z_j^{(R)}) = \p_{w_j} $;
\end{itemize}
see \cite{Folland}.  It is important to recall that, with our choice
of normalization of the Fourier transform, if $f,g\in L^1(\bbH_n)$,
 $\sigma_\lambda(f*g) = 
\sigma_\lambda(f)\sigma_\lambda(g)$, so that 
\begin{equation*}
\sigma_\lambda(V^{(L)} f) 
= - \sigma_\lambda(f) d\sigma_\lambda(V^{(L)})
\quad\text{and}\quad  
\sigma_\lambda(V^{(R)} f) 
= - d\sigma_\lambda(V^{(R)}) \sigma_\lambda(f) \,,
\end{equation*}
where $V^{(L)}$ and $V^{(R)}$ denote a left-invariant and a
right-invariant vector field, respectively.
\ms

If $f\in L^2(\bbH_n)$, we have Plancherel's formula
\begin{equation*}
\| f\|_{L^2(\bbH_n)}^2 = \frac{1}{(2\pi)^{n+1} }
\int_{\bbR}  \|\sigma_\lambda(f)\|_{\HS}^2 |\lambda|^n\,
d\lambda\,,
\end{equation*}
and, if $f\in L^1\cap L^2(\bbH_n)$ the inversion formula
\begin{equation}\label{inve-form}
f[z,t] = \frac{1}{(2\pi)^{n+1}}  
\int_{\bbR}  \tr \big( \sigma_\lambda(f)\sigma_\lambda[z,t]^*\big)
|\lambda|^n\, d\lambda \,.
\end{equation}
\ms

\subsection{The Cauchy--Riemann equations and the Hardy space}\label{Har-subsec} 

We consider now functions that are holomorphic in $\cU$.  
For $F\in\Hol(\cU)$, recalling \eqref{h-coordinates}, we write 
$\widetilde F=F\circ\Psi^{-1}$, so that
$$
F (\z',\z_{n+1}) =
\widetilde F \big( \z',
\textstyle{ \frac{ \z_{n+1}+\ov \z_{n+1}}{2}} ;
\textstyle{ \frac{ \z_{n+1}-\ov \z_{n+1}}{2i} -\frac14 \z' \cdot\bar\z'
} \big)  \,.
$$
The  equation $\p_{\ov \z_{n+1}}F=0$ now reads
\begin{align*}
0 & = \p_{\ov \z_{n+1}}F (\z',\z_{n+1}) =
  \half \p_t   \widetilde F(z,t,h) -\textstyle{\frac{1}{2i}} \p_h 
\widetilde F(z,t,h) \\
& = \half \big( \p_t   \widetilde F +i\p_h   \widetilde F\big)
(z,t,h) \,, 
\end{align*}
that is, 
\begin{equation}\label{partial-h}
i\p_t \widetilde F_h = 
\p_h \widetilde F_h 
\,.
\end{equation}
The remaining Cauchy--Riemann equations 
$\p_{\ov \z_j}F=0$, $j=1,\dots,n$, respectively,
become
\begin{align}
0 & = \p_{\ov \z_j}F (\z',\z_{n+1}) = \big( \p_{\ov z_j}     -\textstyle{\frac i4} z_j \p_t\big) 
\widetilde F_h[z,t] = \ov Z_j^{(R)} \widetilde F_h[z,t] \label{partial-zj}
\,.
\end{align}
If we also have that
 $\widetilde F_h\in L^1(\bbH_n)$, 
using \eqref{partial-zj} 
we obtain that
\begin{align*}
0 & = \sigma_\lambda\big( \ov Z_j^{(R)} \widetilde F_h \big)  
= -d\sigma_\lambda\big( \ov Z_j^{(R)}\big) \sigma_\lambda 
\big( \widetilde F_h\big) \,,
\end{align*}
for $j=1,\dots,n$.
These imply that $\ran \big( \sigma_\lambda 
(\widetilde F_h)\big)\subseteq \ker d\sigma_\lambda\big( \ov Z_j^{(R)}\big)$,
so that, by  (ii) and (iii) in Subsection \ref{Fock-space-subsec}, it follows that
\begin{itemize}
\item[{\tiny$\bullet$}] $\sigma_\lambda 
(\widetilde F_h)=0$ for $\lambda>0$; \smallskip 
\item[{\tiny$\bullet$}] $\ran \big( \sigma_\lambda 
(\widetilde F_h)\big) \subseteq \operatorname{span}\{1\}$. 
\end{itemize}
We learnt this argument from \cite{R}. \ms

We recall that 
$H^2$ defined in \eqref{Hardy-space} is a reproducing kernel Hilbert space, whose inner
 product 
can be realized by the $L^2$-inner
product of the boundary values, that is,
\begin{equation*}
\la F,G\ra_{H^2}
= \int_{\bbH_n} \widetilde F_0 [z,t] \ov{\widetilde G_0[z,t]} \, dzdt \,.  
\end{equation*}
The
reproducing kernel, which
 is called the {\em Szeg\"o} kernel,  is given by
$$
S(\omega,\z) =  \frac{n!}{(4\pi)^{n+1}}  \Big( \frac{\omega_{n+1}-\ov\z_{n+1}}{2i} -
\textstyle{\frac{1}{4}}\omega'\cdot\bar\z' 
\Big)^{-(n+1)} \,,
$$
or, equivalently,
$$
\widetilde S_{(z,t,h),k} [w,s]=
 \frac{n!}{(2\pi)^{n+1}} 
\Big( h+k -i\big(s-t  + \half \Im(w\cdot\ov z) \big)+\textstyle{\frac{1}{4}} |w-z|^2 
 \Big)^{-(n+1)} 
\,,
$$
see e.g. \cite{Stein}.
\ms

\section{The weighted Bergman spaces}\label{Berg-section}

The next result will be used repeatedly throughout the remainder of
the paper. We recall that the spaces $\cL^2_\nu$ where defined in \eqref{L2-nu}.
\begin{lem}\label{NEW}
Let $\nu>-n-2$ and $\tau\in \cL^2_\nu$.  For $(z,t,h)=\Psi(\z)$ with $\z\in\cU$
define
$$
\widetilde F(z,t,h) =
\frac{1}{(2\pi)^{n+1}} \int_{-\infty}^0 e^{\lambda h} \tr
\big(\tau(\lambda) \sigma_\lambda[z,t]^* \big)\, |\lambda|^n d\lambda \,.
$$
Then, $F$ is holomorphic in $\cU$.
\end{lem}

\proof
We first show that the integral defining $F$ converges absolutely.
Consider the orthonormal basis $\{e_\alpha\}$ of $\cF^\lambda$, where
$e_\alpha(z) =z^\alpha/\|z^\alpha\|_{\cF^\lambda}$, $\alpha$ a multiindex.
As $\tau(\lambda)$ and $\sigma_\lambda[z,t]$ are operators on
$\cF^\lambda$, we compute
\begin{align} 
\tr \big(\tau(\lambda) \sigma_\lambda[z,t]^* \big) 
& =\tr \big( \sigma_\lambda[z,t]^* \tau(\lambda)  \big) 
= \sum_\alpha \la \tau(\lambda) e_\alpha ,\,
\sigma_\lambda[z,t] e_\alpha \ra_{\cF^\lambda} \notag \\
& = \sum_\alpha  \la \tau(\lambda) e_\alpha ,\,
P_0 \sigma_\lambda[z,t] e_\alpha \ra_{\cF^\lambda} = \tr \big(\tau(\lambda) P_0\sigma_\lambda[z,t]^* \big) 
\,, \label{traccia}
\end{align}
where $P_0$ denotes the orthogonal projection onto the subspace
generated by $e_0$, 
since $\ran (\tau(\lambda)) \subseteq \operatorname{span}\{e_0\}$.  
Therefore,
$$
\big| \tr (\tau(\lambda) \sigma_\lambda[z,t]^* ) \big|
\le \| \tau(\lambda)\|_{\HS} \| P_0\sigma_\lambda[z,t] 
\|_{\HS} \,.
$$
Since
$\lambda<0$, we have that
\begin{align}
P_0 \sigma_\lambda[z,t] e_\alpha 
& = \la \sigma_\lambda[z,t] e_\alpha , e_0\ra_{\cF^\lambda} e_0  = 
\bigg(
e^{i\lambda t+ \frac\lambda4 |z^2|} 
\bigg( \frac{|\lambda|}{2\pi} \bigg)^n \int_{\bbC^n} 
\frac{e^{\frac\lambda2 w\cdot z} (\ov z+w)^\alpha}{\|w^\alpha\|_{\cF^\lambda}} \,
e^{-\frac{|\lambda|}{2} |w|^2}dw \bigg) e_0 \notag\\
& =  \bigg( \frac{1}{\sqrt{\alpha!}}
\bigg(\frac{|\lambda|}{2}\bigg)^{|\alpha|/2} 
e^{i\lambda t+ \frac\lambda4
  |z|^2} 
\bar z^\alpha \bigg) e_0\,. \label{P0-sigma}
\end{align}
Therefore,
\begin{align*}
\| P_0 \sigma_\lambda[z,t]  \|_{\HS}^2
& = e^{\frac\lambda2  |z|^2} \sum_\alpha  \frac{1}{\alpha!}
\bigg(\frac{|\lambda|}{2}\bigg)^{|\alpha|} 
| z^\alpha|^2 =1 \,,
\end{align*}
so that
\begin{align}
\int_{-\infty}^0 e^{\lambda h} \big|\tr
(\tau(\lambda) \sigma_\lambda[z,t]^* )\big|\, |\lambda|^n d\lambda 
& \le \int_{-\infty}^0 e^{\lambda h} \|\tau(\lambda)\|_{\HS} \,
|\lambda|^n d\lambda \notag \\
& \le \| \tau\|_{\cL^2_\nu} \bigg(\int_{-\infty}^0 e^{2\lambda h}
|\lambda|^{n+\nu+1} d\lambda \bigg)^{1/2} \,, \label{NEW-equation}
\end{align}
which is finite since $\nu>-n-2$. This inequality also shows that the
integral is locally uniformly bounded in $(z,t,h)\in
\bU$. 

In order to show that $F$ is
holomorphic in $\cU$, by the previous estimate, it suffices to show that the
integrand $\widetilde J(z,t,h)=e^{h\lambda} \tr \big(
\tau(\lambda)\sigma_\lambda[z,t]^*\big)$
satisfies equations
\eqref{partial-h} and \eqref{partial-zj}. Indeed, using
\eqref{traccia} and
\eqref{P0-sigma} we have 
\begin{align*}
\widetilde J(z,t,h) 
& =
 e^{h\lambda} \tr \big(
\tau(\lambda) \sigma_\lambda[z,t]^*\big)  \\
& = e^{h\lambda-i\lambda t+ \frac\lambda4
  |z|^2} 
\sum_\alpha  \textstyle{ \frac{1}{\sqrt{\alpha!}}
\big(\frac{|\lambda|}{2}\big)^{|\alpha|/2} }
 z^\alpha 
\la \tau(\lambda) e_\alpha ,\,
 e_0 \ra_{\cF^\lambda}  \,.
\end{align*}
Hence,
$$
(\p_t+i\p_h) \widetilde J (z,t,h) =
(\p_{\ov z_j} -\textstyle{ \frac{i}{4}}z_j \p_t) \widetilde J (z,t,h) 
= 0\,. 
$$
The conclusion follows.
\ms\qed

We now turn to the Bergman spaces.  We begin with the elementary observation that
if $\Phi_{[w,s]}, \Phi_U$ and $D_\delta$ are
  as in Lemma \ref{Aut-U}, and $F\in A^2_\nu$, then
$F\circ \Phi_{[w,s]}$ and $F\circ \Phi_U$ have the same norm as $F$,
while
$\|F\circ D_\delta\|_{A^2_\nu} = \delta^{-(2n+4)/2} \|F\|_{A^2_\nu}
$.  We set $F_{(\eps)} = F(\cdot+ \eps\i)$.

 \begin{prop}\label{prel-prop-Berg-sp}
Let $\nu>-1$.  The following properties hold.
\begin{itemize}
\item[(i)] There exists a constant   $C > 0$ such that
for all $\z \in \cU$, $\eps>0$ and $F \in A_\nu^2 $,
$$
|F(\z+\eps\i)| \leq C\eps^{-(n+2+\nu)/2} \|F\|_{A_\nu^2} \,.
$$ 
As a consequence,  $A^2_\nu$ is a reproducing kernel Hilbert space. 
\item[(ii)]
 There exists a constant  $C> 0$ such that for all
 $\eps>0$ 
and $F \in A_\nu^2 $
 $$ 
\|\widetilde F_{(\eps),h}\|_{L^2(\bbH_n)} 
\le C \eps^{-(1+\nu)/2} \|F\|_{A_\nu^2} \,.
$$ 
In particular, $F_{(\eps)}\in H^2$ and 
$\|F_{(\eps)}\|_{H^2} \le C \eps^{-(1+\nu)/2} \|F\|_{A_\nu^2} $. 
\end{itemize}
\end{prop}

\proof
We begin by observing that if $\z = h\i$, then
 $P( h\i, r)$ is comparable to $P=P(r)=\{ (w,s,k) \; :
|w|<r, \, |s|<r^2,\, |h-k|<r^2\}$. 
We have that
\begin{align*}
|F(h\i)|^2 
& \le C\frac{1}{|P(r)|}\int_{P}|F(\omega)|^2\, d\omega\\
&\le C h^{-(n+2)} \int_{|h-k|<r^2} \int_{B([0,0],r)} 
| \widetilde F (w,s,k)|^2\, dwds \, dk \,.
\end{align*}
For a generic $\z=\Psi^{-1}(z,t,h)$, we have $\z=\Phi_{[z,t]}(h\i)$,
and
$F(\z)= \big(F\circ\Phi_{[z,t]})(h\i)$,
so that
\begin{align}
|F(\z)|^2 
&\le C h^{-(n+2)} \int_{|h-k|<r^2} \int_{B([0,0],r)}
| (F\circ\Phi_{[z,t]})\widetilde{\ } (w,s,k)|^2\, dwds \, dk\notag \\
& =  C h^{-(n+2)} \int_{|h-k|<r^2} \int_{B([z,t]),r)}
| \widetilde F (w,s,k)|^2\, dwds \, dk \notag \\
& =  C h^{-(n+2)} \int_{P((z,t,h),r )}
| \widetilde F (w,s,k)|^2\, dwds \, dk \label{pre-i}
\,.
\end{align}

Given $\eps>0$, we apply \eqref{pre-i} to
$F(\z+\eps\i)=\widetilde F_{(h)}(z,t,\eps)$ and obtain
\begin{align*}
|F(\z+\eps\i)|^2
& = |\widetilde F_{(h)}(z,t,\eps)|^2
\le C \eps^{-(n+2)} \int_{P((z,t,\eps),\sqrt{\eps/2} )}
| \widetilde F_{(h)} (w,s,k)|^2\, dwds \, dk \\
& \le C \eps^{-(n+2+\nu)} \int_{P((z,t,\eps),\sqrt{\eps/2} )}
| \widetilde F_{(h)} (w,s,k)|^2\, dwds \, k^\nu dk \,.
\end{align*}
Since $\| F_{(h)}\|_{A^2_\nu} \le \| F\|_{A^2_\nu} $, this  proves
(i).  
\ms

Next, if $F \in A_\nu^2$ 
and $\eps>0$,
\begin{align*}
|\widetilde F_{h+\eps}[z,t]|^2
 &  = |F(\z+\eps\i)|^2
\le C \eps^{-(n+2+\nu)} \int_{P((z,t,\eps),\sqrt{\eps/2} )}
| \widetilde F_{(h)} (w,s,k)|^2\, dwds \, k^\nu dk  \\
&\le C \eps^{-(n+2+\nu)}\int_{\eps/2}^{3\eps/2} 
\int_{B([z,t],\sqrt{\eps/2})} |\widetilde F_{h+k}[w,s]|^2 \, dwds 
\,k^\nu dk
\,.
\end{align*}
Therefore,
\begin{align*}
\|\widetilde F_{(\eps),h}\|_{L^2(\bbH_n)}^2
&\le C \eps^{-(n+2+\nu)} \int_{\bbH_n}  \int_{\eps/2}^{3\eps/2} 
\int_{B([z,t],\sqrt{\eps/2})} 
|\widetilde F_{h+k}[w,s]|^2 \, dwds \, k^\nu dk  
\,dzdt\\
& = C \eps^{-(n+2+\nu)}  
\int_{B([0,0],\eps)} \int_{\eps/2}^{3\eps/2} 
 \int_{\bbH_n} |\widetilde F_{h+k}\big([w,s][z,t]\big)|^2 \, dzdt \,k^\nu dk
\,dwds \\
& \le C \eps^{-(1+\nu)}  \|F\|_{A^2_\nu}^2 
\,.
\end{align*}
  Moreover,
$$
\| \widetilde F_\eps\|_{L^2(\bbH_n)}\le
\|F_{(\eps)}\|_{H^2}
= \sup_{h>0}  \|\widetilde F_{\eps+h}\|_{L^2(\bbH_n)}
\le C \eps^{-2(1+\nu)/p}  \|F\|_{A^2_\nu} \,.
$$
This proves (ii).
\qed\ms

With a standard argument, see, for instance, \cite[Remark 1.16]{camerun}, it is possible to prove the following result.
\begin{prop}
If $\nu\le -1$, then
$A^2_\nu=\{0\}$.
\end{prop}

We now have all the ingredients to prove our first main result.

\proof[Proof of Theorem \ref{PW-thm-Berg}]
Let $F\in A^2_\nu$ and $\eps>0$.  Then, $F_{(\eps)}\in H^2$ so that
$\widetilde F_{(\eps),h} \in L^2(\bbH_n)$ for $h>0$, and
$\sigma_\lambda (\widetilde F_{\eps+h}) 
=\sigma_\lambda (\widetilde F_{(\eps),h}) =0$ if $\lambda>0$.
Moreover, by \eqref{ext-form} and the inversion formula
\eqref{inve-form}, we have that 
$\sigma_\lambda (\widetilde F_{(\eps),h}) = e^{h\lambda} (\widetilde
F_{(\eps),0})$ if $\lambda<0$. 
Therefore, if $F\in A^2_\nu$, $\eps,h>0$:
\begin{itemize}
\item[{\tiny$\bullet$}] $\widetilde F_h \in L^2(\bbH_n)$ for all
  $h>0$; \smallskip
\item[{\tiny$\bullet$}] $\sigma_\lambda (\widetilde F_h) =0$ if
  $\lambda>0$;\smallskip 
\item[{\tiny$\bullet$}] $\ran( 
\sigma_\lambda (\widetilde F_h) )\subseteq \operatorname{span}\{1\}$ if $\lambda<0$; 
\smallskip
\item[{\tiny$\bullet$}] 
$\sigma_\lambda (\widetilde F_{\eps+h}) =
  e^{h\lambda}
\sigma (\widetilde
F_{\eps})$. 
\end{itemize}

Since $F_{(\eps)}\in H^2$, by the Ogden--Vagi Theorem, there exists
$g_\eps\in L^2(\bbH_n)$ such that
$\sigma_\lambda (g_\eps) =0$ if
  $\lambda>0$,  $\ran( 
\sigma_\lambda (g_\eps) )\subseteq \operatorname{span}\{1\}$ if
$\lambda<0$, and 
\begin{align*}
F(\z+\eps\i) 
& =  \widetilde F_{(\eps), h} [z,t] 
 = \frac{1}{(2\pi)^{n+1}} \int_{-\infty}^0 e^{h\lambda} \tr
\big(\sigma_\lambda(g_\eps) \sigma_\lambda[z,t]^*\big) \, |\lambda|^n
d\lambda \,,
\end{align*}
where $\Psi(\z)=(z,t,h)$.
Switching the roles of $h$ and $\eps$ and
arguing as above, there exists $g_h\in L^2(\bbH_n)$  such that
$\sigma_\lambda (g_h) =0$ if
  $\lambda>0$,  $\ran( 
\sigma_\lambda (g_h) )\subseteq \operatorname{span}\{1\}$ if
$\lambda<0$, and 
$$
F(\z+\eps\i) = 
\frac{1}{(2\pi)^{n+1}} \int_{-\infty}^0 e^{\eps\lambda} \tr
\big(\sigma_\lambda(g_h) \sigma_\lambda[z,t]^*\big) \, |\lambda|^n
d\lambda \,.
$$
These equalities imply that $e^{h\lambda} \sigma_\lambda(g_\eps)  
=  e^{\eps\lambda} \sigma_\lambda(g_h) $ for all $\eps,h>0$, that is,
for every $\lambda<0$
$$
\HS(\cF^\lambda)\ni e^{-\eps\lambda} \sigma_\lambda(g_\eps) =:
\tau(\lambda) 
$$
is well defined, i.e. independent of $\eps$, with $\tau(\lambda)=0$ if
$\lambda>0$ and $\ran(\tau(\lambda))\subseteq
\operatorname{span}\{1\}$.  
Hence, 
$$
F(\z+\eps\i)
= \frac{1}{(2\pi)^{n+1}} 
\int_{-\infty}^0 e^{(h+\eps)\lambda} \tr
\big( \tau(\lambda) \sigma_\lambda[z,t]^*\big) \, |\lambda|^n
d\lambda \,,
$$
that is,
$$
F(\z)
= \frac{1}{(2\pi)^{n+1}} 
 \int_{-\infty}^0 e^{h\lambda} \tr
\big( \tau(\lambda) \sigma_\lambda[z,t]^*\big) \, |\lambda|^n
d\lambda \,.
$$
In particular,
$\sigma_\lambda(\widetilde F_h) = e^{h\lambda}\tau(\lambda)$.
Therefore,
\begin{align}
\| F \|_{A^2_\nu}^2 
& = \int_0^{+\infty} \int_{\bbH_n} |\widetilde F_h [z,t]|^2\,
dzdt \, h^\nu dh \notag \\
& = \frac{1}{(2\pi)^{n+1} } \int_0^{+\infty} 
\int_{-\infty}^0   
 \| \sigma_\lambda(\widetilde F_h) \|_{\HS}^2 \, |\lambda|^n d\lambda
\, h^\nu dh
\notag \\
&
 = \frac{1}{(2\pi)^{n+1} } \int_0^{+\infty} \| \tau(-\lambda)\|_{\HS}^2  
\int_0^{+\infty} e^{-2\lambda h}h^\nu dh\,
\lambda^n d\lambda \notag \\
& = \frac{\Gamma(\nu+1)}{2^{\nu+1}(2\pi)^{n+1}} 
\int_0^{+\infty} \| \tau(-\lambda)\|_{\HS}^2   \lambda^{n-1-\nu}  \,
d\lambda \notag \\
& = \frac{\Gamma(\nu+1)}{2^{\nu+1}} \|\tau\|_{\cL^2_\nu}^2  
\,. \label{Planc-PW-A2nu}
\end{align}

Conversely, 
let $\tau\in\cL^2_\nu$ and $F$ defined by \eqref{PW-A2nu-eq1}.
By Lemma \ref{NEW} 
we
have that $F\in\Hol(\cU)$.  
Moreover, $\widetilde F_h\in L^2(\bbH_n)$ for every $h>0$, since
$\nu>-1$ and by
Plancherel's formula
$$
\| \widetilde F_h\|_{L^2(\bbH_n)}^2 = 
\frac{1}{(2\pi)^{n+1}} \int_{-\infty}^0 \| e^{\lambda h}
\tau(\lambda)\|_{\HS}^2\, |\lambda|^n d\lambda \le C_h
\|\tau\|_{\cL^2_\nu}^2 \,.
$$ 
Moreover, 
$\sigma_\lambda(\widetilde F_h) = e^{h\lambda}\tau(\lambda)$. Hence,
identities \eqref{Planc-PW-A2nu} hold true, and 
\eqref{PW-A2nu-eq2} follows.
\qed
\ms

An immediate consequence is the following result.

\begin{cor}\label{density-Bergman}
Let $\nu>-1$ and $F\in A^2_\nu$.  For $\eps>0$, let $F_{(\eps)}(\z) = F(\z+\eps\i)$.  Then, $F$ is
holomorphic in a neighborhood of $\cU$ and 
$F_{(\eps)}\to F$ in $A^2_\nu$ as $\eps\to0^+$.
\end{cor}

\ms

\begin{remark}\label{weighted-Berg-kernel}{\rm
It is well known that
the reproducing kernel for $A^2_\nu$ is the kernel function, called 
the {\em weighted Bergman  kernel},
$$
K_\nu(\omega,\z)=
\gamma_{n,\nu} \Big( \frac{\omega_{n+1}-\ov\z_{n+1}}{2i} -
\textstyle{\frac{1}{4}}\omega'\cdot\ov{\z'} 
\Big)^{-(n+2+\nu)} \,, 
$$
where
$\gamma_{n,\nu}=\frac{1}{(4\pi)^{n+1}}\frac{\Gamma(n+2+\nu)}{\Gamma(\nu+1)}$.
This fact can be obtained from the expression of the kernel of 
the corresponding weighted Bergman space on the unit ball, by means of the
transformation rule for the Bergman kernel, or as a corollary of the
Paley--Wiener theorem, using the same techniques we will use in Corollary
\ref{repr-ker-weighted-Dir-cor}. 
} \ms
\end{remark}

\section{The weighted Dirichlet spaces}\label{w-Dir-sec}
\ms

Recall that the weighted Dirichlet spaces
$\cD_{\nu,(m)}$ are defined in \eqref{w-Dir-space-def}.
Note that condition (i) in \eqref{w-Dir-space-def} means that
$$
\lim_{\Im \z_{n+1}\to+\infty}\sup_{|\z'|\le R}|F(\z)|=0 \,.
$$
An analogous remark holds for (i) in \eqref{Dir-space-def} as well. 
We begin with an elementary lemma.

\begin{lem}\label{wDs-L1}
The following properties hold true.
\begin{itemize}
\item[(i)] Let $a,b\in\bbR$. If $a>-1$ and $b>0$, then there exists
  $C_0>0$ such that
 $$
I(\z)=\int_\cU \frac{\rho^a(\omega)}{\big|
  \frac{\z_{n+1}-\bar\omega_{n+1} }{2i} - \frac14 \z'\cdot\bar\omega'
    \big|^{a+b+n+2} } \, d\omega = C_0 \frac{1}{(\Im \z_{n+1}-\frac14|\z'|^2)^b} \,.
$$
If  $a\le-1$ or $b\le0$, then the above integral equals $+\infty$.\smallskip
\item[(ii)]   The spaces $\cD_{\nu,(m)}$ are reproducing kernel Hilbert spaces.
\end{itemize}
\end{lem}

\proof (i)  This is an elementary calculation.  We provide
the details for sake of completeness.  We observe that, if
$(z,t,h)=\Psi(\z',\z_{n+1})$ and $(w,s,k)=\Psi(\omega',\omega_{n+1})$, we have
\begin{align*}
\big|
  \textstyle{\frac{1}{2i} }\big(\z_{n+1}-\bar\omega_{n+1} \big) 
  -  \textstyle{\frac14} \z'\cdot\bar\omega'
    \big|^2  
& =   \textstyle{\frac14 }
\big(
  h+k +\textstyle{\frac12} |z- w|^2 \big)^2 
 + \textstyle{\frac14 } \big(  (s-t) 
-  \textstyle{\frac{1}{2} } \Im z\cdot\bar w
    \big)^2  
\,.
\end{align*}
 Then, by the standard translation invariance of the Lebesgue measure
 in $\bbR$ and $\bbC^n$, and integration in polar coordinates in
 $\bbC^n$,  we have
\begin{align*}
\widetilde I(z,t,h) & = 2^{a+b+n+2}  \int_\bU \frac{k^a}{
  \big(\big(h+k+\frac14|w|^2\big)^2+s^2\big)^{(a+b+n+2)/2} } \, dsdwdk
\\
& = \frac{2^{a+b+n+3} \pi^{n+1}}{n!} 
\int_0^{+\infty} \int_0^{+\infty}  \int_\bbR \frac{k^a r^{2n-1}}{
\big(\big(h+k+\frac14 r^2\big)^2+s^2\big)^{(a+b+n+2)/2} } \, dsdrdk \\
& = \frac{2^{a+b+n+3} \pi^{n+1}}{n!} C_1
\int_0^{+\infty} \int_0^{+\infty}  \frac{k^a r^{2n-1}}{
\big(h+k+\frac14 r^2\big)^{a+b+n+1 }} \, drdk \\
& = \frac{2^{a+b+3n+3} \pi^{n+1}}{n!} C_1 C_2 \int_0^{+\infty}  \frac{k^a}{
\big(h+k\big)^{a+b+1}} \, dk \\
& = \frac{C_0}{h^b} \,,
\end{align*}
as we wished to show.

(ii) Let $F\in\cD_{\nu,(m)}$.  Then, $\p_{\z_{n+1}}^m F\in
A^2_{2m+\nu}$.  For $\z\in\cU$ we define
\begin{equation} \label{G-def}
G(\z)
 = c \int_\cU \frac{\p_{\z_{n+1}}^m F(\omega)}{\big(
  \frac{\z_{n+1}-\bar\omega_{n+1} }{2i} - \frac14 \z'\cdot\bar\omega'
    \big)^{n+2+\nu+m} }\,\rho(\omega)^{2m+\nu}  d\omega \,,
\end{equation} 
where
$c$ is a constant to be chosen later.
Then, (i) and Cauchy--Schwarz's inequality give that
\begin{equation} \label{G-est}
|G(\z)|  \le C_0 h^{-(n+2+\nu)/2}  
\big\| \p_{\z_{n+1}}^m F\|_{A^2_{2m+\nu}} \,, 
\end{equation} 
where $C_0$ is as in (i), and $h=\Im \z_{n+1}-\frac14 |\z'|^2$.
Arguing as above, it is easy to see that $G$ is holomorphic in 
$\cU$, and that we can differentiate under the integral sign $m$ times to obtain
that, using  Remark
\ref{weighted-Berg-kernel}, for a suitable constant $c$,  
$\p_{\z_{n+1}}^m G=\p_{\z_{n+1}}^m F$.  Therefore, 
$(F-G)(\z)=\sum_{j=0}^{m-1} g_j(\z') \z_{n+1}^j$, where the
$g_j$'s are entire 
functions in $\bbC^n$.  
By \eqref{G-est} it also follows that 
$\lim_{|\z'|\le R,\, \Im\z_{n+1}\to+\infty} G(\z) =0$.  Therefore,
for each $\z'$ fixed, the polynomial $F(\z',\cdot)-G(\z',\cdot)$ tends
to $0$ as $\Im\z_{n+1}\to+\infty$.  This implies that
the $g_j$'s are identically $0$; hence $G=F$. 

Now, \eqref{G-est} with $G$ replaced by $F$ gives that the point
evaluations are bounded on $\cD_{\nu,(m)}$, and also implies
 uniform estimates on
compact subsets of $\cU$. An elementary argument 
shows that $\cD_{\nu,(m)}$
is complete; hence a reproducing kernel Hilbert space.  
\ms\qed

We set
\begin{equation} \label{H-m-def}
\cH_m
=\big\{  F \in\Hol(\ov\cU):\, 
\p_\z^\alpha
F \in H^2,\ |\alpha|\le m\big\}\,.
\end{equation}

\begin{lem}\label{wDs-density-lem}
Let $-(n+2)< \nu<-1$, and let $m>-\frac{\nu+1}{2}$.
  Then, $\cD_{\nu,(m)}\cap \cH_m$ is dense in
$\cD_{\nu,(m)}$.
\end{lem}

\proof
Let $F\in\cD_{\nu,(m)}$.  
For $\eps,\delta>0$, $q>0$ to be selected later, and $\z\in\cU$, we
define  
\begin{equation*} 
G_{(\eps,\delta)} (\z)
 = c \int_\cU 
\frac{1}{\big( -\eps i\omega_{n+1} +1\big)^q} 
\cdot\frac{\p_{\z_{n+1}}^m F(\omega)}{\big(
  \frac{\z_{n+1}+i\delta- \bar\omega_{n+1}  }{2i} - \frac14 \z'\cdot\bar\omega'
    \big)^{n+2+\nu+m} }\,\rho(\omega)^{2m+\nu}  d\omega \,,
\end{equation*} 
where $c$ is as in \eqref{G-def}.
Recall that $\p_{\z_{n+1}}^m F\in
A^2_{2m+\nu}$.
Observe that the factor 
$\big( -\eps i\omega_{n+1} +1\big)^{-q} $
is bounded on $\ov\cU$.  Then, the same
argument as in Lemma \ref{wDs-L1} (ii) and  the dominated convergence
theorem  give that  $G_{(\eps,\delta)}$ are holomorphic and, by
\eqref{G-def}, 
converge
uniformly on compact subsets to $F$, as $\eps,\delta\to0^+$. 
 Thus, we need to show: (a)  that
$G_{(\eps,\delta)} \in\cH_m$; and (b) that converge to $F$ in $\cD_{\nu,(m)}$.  
\ms

Let $\alpha=(\alpha',\alpha_{n+1})$ be a multiindex, $|\alpha|\le m$.  Then,
differentiating under the integral sign, for a suitable constant $c'$, 
we have
\begin{align*}
\p_\z^\alpha G_{(\eps,\delta)} (\z)
& =c' 
\int_\cU 
\frac{{\ov{\omega}'}^{\alpha'}}{\big( -\eps i\omega_{n+1} +1\big)^q} 
\cdot\frac{\p_{\z_{n+1}}^m F(\omega)}{\big(
  \frac{\z_{n+1}+i\delta-\bar\omega_{n+1} }{2i} - \frac14 \z'\cdot\bar\omega'
    \big)^{n+2+\nu+m+|\alpha|} }\,\rho(\omega)^{2m+\nu}\,  d\omega 
 \\
& =: \int_\cU K(\z,\omega) \p_{\z_{n+1}}^m F(\omega)
\rho(\omega)^{2m+\nu}\,  d\omega 
\,.
\end{align*}  

Letting 
$(z,t,h)=\Psi(\z',\z_{n+1})$, $(w,s,k)=\Psi(\omega',\omega_{n+1})$ and
writing $\widetilde K = K \big(\Psi(\cdot),\Psi(\cdot)\big)$, we see
that
\begin{align*}
& |\widetilde K \big( (z,t,h),(w,s,k)\big)|  \\
& \le C
\frac{1}{ \big( \eps (k+\frac14|w|^2 +|s|)+ 1 \big)^{q-|\alpha'|/2} } \cdot
\frac{1}{  
\big(
  h+\delta+ k +\textstyle{\frac12} |z- w|^2 
 + \big|  (s-t) 
-  \textstyle{\frac{1}{2} } \Im z\cdot\bar w
    \big| \big)^{n+2+\nu+m+|\alpha|} } \\
& =: C \widetilde L \big( (z,t,h),(w,s,k)\big)
 \,.
\end{align*}

Using Cauchy--Schwarz's inequality, for $q>0$ sufficiently large, we have
\begin{align}
 & \int_{\bbH_n} \big| \big(\p_\z^\alpha G_{(\eps,\delta)}\big)_h 
\!\!\widetilde{\,\,}\, [z,t]\big|^2\, dzdt \notag\\
& \le \int_{\bbH_n}  \bigg(
\int_\bU \widetilde L \big( (z,t,h),(w,s,k)\big)
|\p_k^m \widetilde F(w,s,k)| 
\, dsdw\, k^{2m+\nu}dk 
\bigg)^2 \, dzdt \notag\\
& \le  \bigg(
\int_\bU 
\frac{1}{ \big( \eps (k+\frac14|w|^2+|s|)+ 1 \big)^{2q-|\alpha'|} } 
, dsdw\, k^{2m+\nu}dk 
\bigg) \notag\\
& \qquad \times
\int_{\bbH_n} 
\int_\bU
\frac{|\p_k^m \widetilde F(w,s,k)|^2}{  
\big(
  h+\delta+ k +\textstyle{\frac12} |z- w|^2 
 + \big|  (s-t) 
-  \textstyle{\frac{1}{2} } \Im z\cdot\bar w
    \big| \big)^{2(n+2+\nu+m+|\alpha|)}  } \, dsdw\, k^{2m+\nu}dk 
dzdt \notag\\
& \le C \int_\bU
|\p_k^m \widetilde F(w,s,k)|^2 
\int_{\bbH_n} \frac{1}{  
\big(
  h+\delta+ k +\textstyle{\frac12} |z|^2 
 + |t| \big)^{2(n+2+\nu+m+|\alpha|)}  } 
\, dzdt  dsdw\, k^{2m+\nu}dk\,.  \label{intermediate}
\end{align} 
By Lemma \ref{wDs-L1} (i) it follows that 
$
\| \p_\z^\alpha G_{(\eps,\delta)} \|_{ H^2}
 \le C \| F\|_{\cD_{\nu,(m)}}^2$, 
for $|\alpha|\le m$,  and  when $|\alpha|=m$, also that
\begin{multline*}
\int_0^{+\infty} \int_{\bbH_n} \big| \big(\p_\z^\alpha G_{(\eps,\delta)}\big)_h 
\!\!\widetilde{\,\,}\, [z,t]\big|^2\, dzdt\, h^{2\nu+m}dh \\ 
\le
C \int_\bU
|\p_k^m \widetilde F(w,s,k)|^2 \bigg( 
\int_\bU \frac{1}{  
\big(
  h+\delta+ k +\textstyle{\frac12} |z|^2 
 + |t| \big)^{2(n+2+\nu+2m)}  } 
\, dzdt \, h^{2\nu+m}dh\bigg)  dsdw\, k^{2m+\nu}dk \,.
\end{multline*}
This implies that $\p_{\z_{n+1}}^m G_{(\eps,\delta)} \in A^2_{2m+\nu}$.
It is  also easy to see that $\lim_{|z|\le R,\, h\to+\infty} \widetilde G_{(\eps,\delta)}(z,t,h)=0$.
Therefore, $G_{(\eps,\delta)} \in \cH_m$, i.e.
the conclusion (a) follows.  Now, it is elementary to show (b).
Indeed, since $\big( -\eps i\omega_{n+1} +1\big)^{-q} 
\p_{\z_{n+1}} ^m F\in A^2_{2m+\nu}$, we have that
\begin{align*}
\p_{\z_{n+1}}^m G_{(\eps,\delta)} (\z) 
& =   c \int_\cU 
\frac{\p_{\z_{n+1}}^m F(\omega)}{\big( -\eps i\omega_{n+1}
  +1\big)^q}
\cdot\frac{1}{\big(
  \frac{\z_{n+1}+i\delta- \bar\omega_{n+1}  }{2i} - \frac14 \z'\cdot\bar\omega'
    \big)^{n+2+2m+\nu} }\,\rho(\omega)^{2m+\nu}  d\omega \\
& = \frac{\p_{\z_{n+1}}^m F(\z+\delta\i)}{\big( -\eps i\z_{n+1} +\eps\delta
  +1\big)^q}\,.
\end{align*}
A simple application of the dominated convergence theorem together
with Corollary \ref{density-Bergman} give the
desired conclusion.
\ms\qed

\proof[Proof of Theorem \ref{PW-thm-Dir}]
 Let $F\in \cD_{\nu,(m)}\cap \cH_m$, with  $m>-\frac{\nu+1}{2}$.  
Observe that, since $\p_{\z_{n+1}}^j F\in A^2_{2m+\nu}\cap H^2$ for
$j=0,\dots,m$, we have that 
$\p_h^j \widetilde F_h \in L^2(\bbH_n)$ for all $h>0$, and
$j=0,\dots,m$.  This easily implies that, for $F \in\cH_m$, setting
$\tau(\lambda) = \sigma_\lambda(\widetilde F_0)$, 
\begin{equation*}
\sigma_\lambda\big( \p_h^m \widetilde F_h\big)
= \p_h^m \sigma_\lambda\big( \widetilde F_h\big)
=  \p_h^m \big(e^{h\lambda} \tau(\lambda) \big)
= \lambda^m e^{h\lambda} \tau(\lambda)
\,.
\end{equation*}
Observe that with this choice of $\tau$, formula \eqref{PW-Dnu-eq1}
holds for $F$.  Moreover,
\begin{align*}
\|F\|_{\cD_{\nu,(m)}}^2
& = \int_\cU |\rho^m(\z) \p_{\z_{n+1}}^m F(\z)|^2 \,
\rho^\nu(\z) d\z \\
& = \int_0^{+\infty} \int_{\bbH_n}  \big|h^m \p_h^m \widetilde F_h[z,t] \big|^2\, dzdt\, h^\nu
dh \\ 
& =  \frac{1}{(2\pi)^{n+1}} 
\int_0^{+\infty} \int_{-\infty}^0   \big\| 
\sigma_\lambda\big( \p_h^m \widetilde F_h\big)
\big\|_{\HS}^2\,|\lambda|^n d\lambda \,  h^{2m+\nu} \, dh \\ 
& = \frac{1}{(2\pi)^{n+1}}  
\int_0^{+\infty}  \int_{-\infty}^0   e^{2h\lambda}
\|\tau(\lambda)\|_{\HS}^2\,|\lambda|^{2m+n} d\lambda \, h^{2m+\nu}\, dh \\ 
& = \frac{1}{(2\pi)^{n+1}} 
\int_0^{+\infty} \|\tau(-\lambda)\|_{\HS}^2 
|\lambda|^{2m+n} \int_0^{+\infty}
e^{-2h\lambda}h^{2m+\nu} \, dh
d\lambda \\ 
& = \frac{1}{(2\pi)^{n+1}} \frac{\Gamma(2m+\nu+1)}{2^{2m+\nu+1}} 
\int_0^{+\infty} \|\tau(-\lambda)\|_{\HS}^2 
|\lambda|^{n-(\nu+1)}  \, 
d\lambda \\
& = \frac{\Gamma(2m+\nu+1)}{2^{2m+\nu+1}} 
\| \tau\|_{\cL^2_\nu}^2
\,.
\end{align*}

Suppose now that $F\in \cD_{\nu,(m)}$ and let $\{F_N\}$ be a sequence
in $\cD_{\nu,(m)}\cap\cH_m$, $F_N\to F$ in $\cD_{\nu,(m)}$.  Then, $F_N$ converges to
$F$ also uniformly on compact subsets.  Let
$\tau_N\in \cL^2_\nu$ be such that 
$ (\widetilde F_N)_h[z,t] =
\frac{1}{(2\pi)^{n+1}} \int_{-\infty}^0 
 e^{h\lambda} \tr \big(
\tau_N(\lambda)\sigma_\lambda[z,t]^*\big) \, |\lambda|^n d\lambda$,
and let $\tau=\lim_{N\to+\infty} \tau_N$ in $\cL^2_\nu$.  
Then,
\begin{align*}
\widetilde F_h[z,t] 
& =  \lim_{N\to+\infty} \frac{1}{(2\pi)^{n+1}} \int_{-\infty}^0 
 e^{h\lambda} \tr \big(
\tau_N(\lambda)\sigma_\lambda[z,t]^*\big) \, |\lambda|^n d\lambda \\
& = \frac{1}{(2\pi)^{n+1}} \int_{-\infty}^0 
 e^{h\lambda} \tr \big(
\tau(\lambda)\sigma_\lambda[z,t]^*\big) \, |\lambda|^n d\lambda \,,
\end{align*}
by applying estimate \eqref{NEW-equation}.
This proves \eqref{PW-Dnu-eq1}, and \eqref{PW-Dnu-eq2} follows as well.

Conversely, let $F$ be given by \eqref{PW-Dnu-eq1}.  Then, Lemma
\ref{NEW} gives that $F$ is holomorphic in $\cU$.  Plancherel's
formula now  gives \eqref{PW-Dnu-eq2}.

Finally, we observe that by \eqref{PW-Dnu-eq2} it follows easily  that
the spaces $\cD_{\nu,(m)}$ do not depend on the choice of the integer
$m$ and their norms coincide, up to a multiplicative constant.
\ms
\qed

\begin{cor}\label{repr-ker-weighted-Dir-cor}
Let $-(n+2)<\nu<-1$.  Let $m$ be a positive integer, $m>-(\nu+1)/2$.  Then,
there exists a constant $\gamma_{\nu,n,m}$ such that  the reproducing
kernel 
$K_\nu$, expressed with respect to the inner product in $\cD_{\nu,(m)}$ is given by
\begin{equation*}
K_\nu(\omega,\z) = \gamma_{n,m,\nu} 
\Big( \frac{\omega_{n+1}-\ov\z_{n+1}}{2i} -
\textstyle{\frac{1}{4}} \omega'\cdot \ov{\z'} 
\Big)^{-(n+2+\nu)} \,,
\end{equation*}
where $\gamma_{n,m,\nu}=
\frac{4^m}{(4\pi)^{n+1}}\frac{\Gamma(n+2+\nu)}{\Gamma(2m+\nu+1)}$.
\end{cor}
\ms

\proof
We use the inversion formula and the polarized identity coming from
the Paley--Wiener type Theorem \ref{PW-thm-Dir}.  For
$F\in\cD_{\nu,(m)}$, let $\tau_F$ denote the element of $\cL_\nu^2$ such
that
$$
\widetilde F(z,t,h) 
 =\frac{1}{(2\pi)^{n+1}} \int_{-\infty}^0 
 e^{h\lambda} \tr \big(
\tau_F(\lambda)P_0\sigma_\lambda[z,t]^*\big) \, |\lambda|^n d\lambda
\,,
$$
where $P_0$ denotes the orthogonal projection onto the subspace
generated by $e_0$.
Moreover, by the reproducing formula for $\cD_{\nu,(m)}$,  
\eqref{PW-Dnu-eq2}, writing
$K_\nu(\z,\cdot)=K_{\z}$ and $(z,t,h)=\Psi(\z)$,
 we have
$$
\widetilde F(z,t,h) = F(\z)
= \la F, K_{\z} \ra_{\cD_{\nu,(m)}} =
\frac{\Gamma(2m+\nu+1)}{2^{2m+\nu+1}(2\pi)^{n+1}} 
\int_{-\infty}^0 \tr \big( \tau_F(\lambda) \tau_{K_\z}(\lambda)^*\big)
\, |\lambda|^{n-\nu-1} d\lambda\,.
$$
Since these two equalities hold for all $\tau\in\cL^2_\nu$, it follows that
$$
\tau_{K_\z}(\lambda) =   \frac{2^{2m+\nu+1}}{\Gamma(2m+\nu+1)} 
|\lambda|^{\nu+1} P_0\sigma_\lambda[z,t]\,.
$$
Therefore, using \eqref{P0-sigma} and writing 
$C=\frac{2^{2m+\nu+1}}{\Gamma(2m+\nu+1)(2\pi)^{n+1}}  $, we have 
\begin{align*}
\widetilde K_{(z,t,h)}(w,s,k)
& = C  \int_{-\infty}^0 
 e^{(h+k)\lambda} \tr \big(
P_0\sigma_\lambda[z,t] P_0 \sigma_\lambda[w,s]^*\big) \, |\lambda|^{n+\nu+1}
d\lambda \\
& = C  \int_{-\infty}^0 
 e^{(h+k)\lambda}
\sum_\alpha \frac{1}{\alpha!}
\bigg(\frac{|\lambda|}{2}\bigg)^{|\alpha|} 
e^{i\lambda (t-s)+ \frac\lambda4
  |z|^2+\frac\lambda4 |w|^2} 
w^\alpha \bar z^\alpha  \, |\lambda|^{n+\nu+1}
d\lambda  \\
& = C  \int_0^{+\infty} 
 e^{-\lambda(h+k -i (s-t+\frac12\Im(w\cdot\bar z )) +\frac14 |w-z|^2) }
\, \lambda^{n+\nu+1}
d\lambda \\
&  = C\, 
\Gamma(n+2+\nu) \Big( h+k +\textstyle{\frac14} |w-z|^2 -i (s-t+ \textstyle{\frac12}\Im(w\cdot\bar z ))
 \Big)^{-(n+2+\nu)} \,,
\end{align*}
that is,
$$
K_\nu(\omega,\z)
= \frac{4^m}{(4\pi)^{n+1}}\frac{\Gamma(n+2+\nu)}{\Gamma(2m+\nu+1)}
\Big( \frac{\omega_{n+1}-\ov\z_{n+1}}{2i} -
\textstyle{\frac{1}{4}} \omega'\cdot \ov{\z'} 
\Big)^{-(n+2+\nu)} \,,
$$
as we wished to show.
\qed
\ms

\section{The Dirichlet space}\label{Dir-sec}
\ms 
In this section we prove Theorem \ref{main-4} and we provide
justification of the name Dirichlet space for the space $\cD_{(m)}$. 

In order to simplify some formulas, we introduce the notation
$$
Q(\omega,\z)=
  \frac{\omega_{n+1}-\ov\z_{n+1}}{2i} - \textstyle{\frac14}
  \omega'\cdot\ov\z'\,,
$$
whereas we remind the reader that $\i$ denotes the point $(0,i)\in\cU$.
\ms

We start proving a couple of lemmas which are the analogue of Lemma
\ref{wDs-L1} and Lemma 
\ref{wDs-density-lem}.
\begin{lem}\label{Ds-L1}
The following properties hold true.
\begin{itemize}
\item[(i)] Let $m>\frac{n+1}{2}$ be fixed. Then, there exists a constant $C>0$ such that
$$
I(\z)=\int_\cU
\bigg|\frac{1}{Q^m(\z,\omega)}-\frac{1}{Q^m(\i,\omega)}\bigg|^2 \rho^{2m-n-2}(\omega)
\, d\omega \leq
C\frac{(1+|\z|^2)^{2m+1}}{(\Im\z_{n+1}-\frac14|\z'|^2)}\,. 
$$
\smallskip
\item[(ii)]   The spaces $\cD_{(m)}$ are reproducing kernel Hilbert spaces.
\end{itemize}
\end{lem}

\proof (i) Given any $\z\in\cU$ there
exists a constant $C>0$, such that 
$$
\big|Q(\i,\omega)-Q(\z,\omega)\big| =
\Big|\frac{i-\z_{n+1}}{2i}
+\frac{1}{4}\z'\cdot\overline{\omega'}\Big|
\le C (1+|\z|)(1+|\omega'|) \,,
$$
so that,
\begin{align*}
\bigg|\frac{1}{Q^m(\z,\omega)}-\frac{1}{Q^m(\i,\omega)}\bigg|^2
& \le C \frac{(1+|\z|)^2(1+|\omega'|)^2}{\big|
  Q^m(\z,\omega)Q^m(\i,\omega)\big|}
\sum_{j=0}^{m-1} 
\big|Q^j(\z,\omega)Q^{m-1-j}(\i,\omega)\big|^2 \,.
\end{align*}
Thus, in order to conclude the proof, 
it is enough to estimate the integral
$$
I_j(\z)=\int_\cU
(1+|\omega'|)^2\bigg|\frac{Q^j(\z,\omega)
Q^{m-1-j}(\i,\omega)}{Q^m(\z,\omega)Q^m(\i,\omega)}\bigg|^2
\rho^{2m-n-2}(\omega)\, d\omega\,. 
$$
Observing that 
$$
\bigg|\frac{Q(\z,\omega)}{Q(\i,\omega)}\bigg| \le C(1+|\z|^2) \,,
$$
it holds that 
\begin{align*}
 I_j(\z)
&=\int_\cU\frac{(1+|\omega'|)^2}{|Q(\i,\omega)|}\bigg|\frac{Q^j(\z,\omega)
   Q^{m-1-j}(\i,\omega)}{Q^m(\z,\omega)Q^{m-\frac12}(\i,\omega)}\bigg|^2
\rho^{2m-n-2}(\omega)\, d\omega\\
 &\leq C(1+|\z|^2)^{1+2j}\int_{\cU}\frac{\rho^{2m-n-2}(\omega)}{|Q(\z,\omega)|^{2m+1}}\, d\omega\\
 &\leq C \frac{(1+|\z|^2)^{1+2j}
 }{(\Im\z_{n+1}-\frac14|\z'|^2)} \,,
\end{align*}
where the last inequality follows from (i) in Lemma \ref{wDs-L1}.
Thus,
\begin{align*}
 I(\z)&\leq C (1+|\z|^2)^2 \sum_{j=0}^{m-1}I_j(\z)
\le C\frac{(1+|\z|^2)^{2m+1}}{(\Im\z_{n+1}-\frac14|\z'|^2)} \,,
\end{align*}
as we wished to prove.
\ms

(ii) Let $F\in \cD_{(m)}$. Then $\p^m_{\z_{n+1}}F\in A^2_{2m-n-2}$. For $\z\in\cU$ we define
\begin{equation}\label{G-def-D}
G(\z)= c\int_\cU \p^m_{\z_{n+1}}F(\omega)\bigg[\frac{1}{Q^m(\z,\omega)}-\frac{1}{Q^m(\i,\omega)}\bigg] \rho(\omega)^{2m-n-2}\, d\omega\,,
\end{equation}
where $c$ is a suitable constant to be chosen later.
Then, (i) and
Cauchy--Schwarz's inequality guarantee that $G$ is well-defined, in
particular, 
\begin{equation}\label{G-D-est}
|G(\zeta)|\leq C\frac{(1+|\z|^2)^{m+\frac12}}{(\Im\z_{n+1}-\frac14|\z'|^2)^{\frac12}}\|  \p_{\z_{n+1}}^m F
\|_{A^2_{2m-n-2}}\leq C\frac{(1+|\z|^2)^{m+\frac12}}{(\Im\z_{n+1}-\frac14|\z'|^2)^{\frac12}}\|  F
\|_{\cD_{(m)}}\,.
\end{equation}
Arguing as in the proof of Lemma \ref{wDs-L1} we obtain that $G\in
\cD_{(m)}$ and, for a suitable choice of the constant $c$, 
$\p^m_{\z_{n+1}}G=\p^m_{\z_{n+1}}F$. Therefore, for each $\z'$ fixed,
we obtain that $(F-G)(\z)= \sum_{j=0}^{m-1}g_j(\z')\z^j_{n+1}$, where
the $g_j$'s are entire functions in $\bbC^n$. Since both
$G$ and $F$ belong to $\cD_{(m)}$, it follows that
$F(\z)-G(\z)= F(\i)$. This fact and
\eqref{G-def-D} give an integral representation for any function $F\in
\cD_{(m)}$ and that \eqref{Dir-norm-def} is a norm.
Finally, this integral representation and \eqref{G-D-est}
show that the point evaluations are bounded on $\cD_{(m)}$, and the fact that
$\cD_{(m)}$ is a reproducing kernel Hilbert space follows as in 
Lemma \ref{wDs-L1}.
\ms 
\qed

\begin{lem}\label{Ds-density-lem}
Let $m>\frac{n+1}{2}$ and let $\cH_m$ be as in \eqref{H-m-def} and let
$\cD_{(m)}(\i)$ be the closed subspace of $\cD_{(m)}$ of functions that
vanish in $\i$. 
Then, $\cD_{(m)}(\i)\cap\cH_m$ is dense in
$\cD_{(m)}(\i)$.
\end{lem}

\proof
Let $F\in \cD_{(m)}(\i)$.  For $\eps,\delta>0$, and $q>0$ to be selected
later, and $\z\in\cU$ we define 
$$
G_{(\eps,\delta)}(\z)=c\int_\cU 
\frac{\p^m_{\z_{n+1}}F(\omega)}{(1-\eps
  i\omega_{n+1})^q}\bigg[\frac{1}{Q^m(\z+\delta\i,\omega)}-\frac{1}{Q^m((1+\delta)\i,\omega)}
\bigg]
\rho(\omega)^{2m-n-2}\, d\omega\,, 
$$
where $c>0$ is as in \eqref{G-def-D}. From Lemma \ref{Ds-L1} and the
dominated convergence theorem we deduce that the functions
$G_{(\eps,\delta)}$ are holomorphic in $\ov\cU$ and 
 converge  to $F$ uniformly on compact subsets, as $\eps,\delta\to0^+$.

We now show that $G_{(\eps,\delta)}\in\cH_m$. Let
$\alpha=(\alpha',\alpha_{n+1})$ be a multiindex,  
$\alpha\leq m$. Then, for a suitable constant $c'$,
\begin{align*}
\p_\z^\alpha G_{(\eps,\delta)} (\z)
& =c' 
\int_\cU 
\frac{{\ov{\omega}'}^{\alpha'}}{\big( -\eps i\omega_{n+1} +1\big)^q} 
\cdot\frac{\p_{\z_{n+1}}^m
  F(\omega)}{Q^{m+|\alpha|}(\z,\omega)}\,\rho(\omega)^{2m-n-2}\,
d\omega  
\notag \\
& =: \int_\cU K(\z,\omega) \p_{\z_{n+1}}^m F(\omega)
\rho(\omega)^{2m-n-2}\, 
\, d\omega 
\,.
\end{align*}  
Arguing as in the proof of Lemma \ref{wDs-density-lem} and selecting
$q$ sufficiently large, we obtain that
$\p^\alpha_\z G_{(\eps,\delta)}\in H^2$ for $|\alpha|\leq m$,
$\p^m_{\z_{n+1}}G_{(\eps,\delta)}\in A^2_{2m-n-2}$ and $\lim_{|\z'|\le
  R,\, \Im \z_{n+1}\to+\infty}\p_{\z_j} G_{(\eps,\delta)}(\z)=0$, for
$j=1,\ldots,n+1$. Thus, $G_{(\eps,\delta)}$ belongs to
$\cD_{(m)}\cap\cH_m$. 
The convergence of $G_{(\eps,\delta)}$ to $F$ in
$\cD_{(m)}$ follows with the same argument as in the
proof of Lemma \ref{wDs-density-lem}. The
proof is therefore complete.
\ms
\qed

In order to prove Theorem \ref{main-4} we need the analogue of Lemma
\ref{NEW}, in the case $\nu=-n-2$.

\begin{lem}\label{NEW2}
Let $\tau\in \cL^2_{-n-2}$.  For $(z,t,h)=\Psi(\z)$ with $\z\in\cU$
define
$$
\widetilde F(z,t,h) =
\frac{1}{(2\pi)^{n+1}} \int_{-\infty}^0
\tr\Big(\tau(\lambda)\big(e^{\lambda
  h}\sigma_\lambda[z,t]^*-e^\lambda\sigma_\lambda[0,0]^*\big) \Big)\,
|\lambda|^n d\lambda\,.
$$
Then, $F$ is holomorphic in $\cU$.
\end{lem}

\proof
Firstly we show that $F$ is defined by an absolutely convergent
integral. From \eqref{traccia} and \eqref{P0-sigma} we obtain that 
\begin{align*}
 \Big|\tr\Big(\tau(\lambda)\big(e^{\lambda h}\sigma_\lambda[z,t]^*-
 &e^\lambda\sigma_\lambda[0,0]^* \big) \Big)\Big|\\ 
 &\leq\Big|\big<\tau(\lambda)e_0,(e^{\lambda(h+it+\frac{|z|^2}{4})
 }-e^\lambda)e_0\big>_{\cF^\lambda}+\sum_{\alpha\neq0}\la
 \tau(\lambda) e_\alpha ,\, 
P_0 \sigma_\lambda[z,t] e_\alpha \ra_{\cF^\lambda} \Big|\\
&\leq
\|\tau(\lambda)\|_{\HS} \Big(\big|e^{\lambda(h+it+\frac{|z|^2}{4})
}-e^\lambda\big|+\big(1-e^{\frac{\lambda}{2}|z|^2}\big)^{\frac12}\Big) \,.
\end{align*}
Therefore, using \eqref{NEW-equation} we have
\begin{align*}
&  \int_{-\infty}^0 
\tr\Big(\tau(\lambda)\big(e^{\lambda
  h}\sigma_\lambda[z,t]^*-e^\lambda\sigma_\lambda[0,0]^*\big) \Big)
\, |\lambda|^n d\lambda\\
 &\leq \| \tau\|_{\cL^2_{-n-2}} \bigg(\int_{-1}^0
 \Big(\big|e^{\lambda(h+it+\frac{|z|^2}{4})
 }-e^\lambda\big|+\big(1-e^{\frac{\lambda}{2}|z|^2}\big)^{\frac12}\Big)^2\,
 |\lambda|^{-1}d\lambda \bigg)^{1/2} \\ 
& \qquad \qquad 
+\|\tau\|_{\cL^2_{-n-2}}\bigg(\int_{-\infty}^{-1}
\big(e^{\lambda h} + e^{\lambda} \big)^2 \,
|\lambda|^{-1}d\lambda \bigg)^{1/2}\\
& <\infty.
\end{align*}
These two last inequalities also show that $\widetilde F(z,t,h)$ is
locally uniformly bounded in $(z,t,h)\in 
\bU$.  
The holomorphicity of $F$ follows arguing as in the proof of Lemma
\ref{NEW}. \ms
\qed

We can now prove Theorem
\ref{main-4}. 
\proof[Proof of Theorem \ref{main-4}.]
We first assume that $F\in \cD_{(m)}(\i)$. Then, from Lemma
\ref{Ds-density-lem} and a minor modification of the proof of Theorem
\ref{PW-thm-Dir}, we obtain \eqref{PW-D-eq1} and \eqref{PW-D-eq2}. If
$F\in \cD_{(m)}$ does not vanish in $\i$, we apply the proof to
$F-F(\i)$ and we are done. 

Conversely, let $F$ be given by \eqref{PW-D-eq1}. Then, Lemma
\ref{NEW2} guarantees the holomorphicity of $F$ in $\cU$ and
Plancherel's formula gives \eqref{PW-D-eq2}. 
\ms
\qed

\begin{cor}\label{repr-ker--Dir-cor}
Let $m$ be a positive integer, $m>(n+1)/2$.  Then, the reproducing
kernel 
$K$, expressed with respect to the inner product in $\cD_{(m)}$ is given by 
\begin{equation*}
K(\omega,\z)= 1+\gamma_{n,m}\log\frac{Q(\omega,\i)Q(\i, \z)}{Q(\omega,\z)}\,,
\end{equation*}
where $\gamma_{n,m}=\frac{2^{2m-n}}{\Gamma(2m-n-1)(2\pi)^{n+1}}$.
\end{cor}
\proof
Let $\tau_F$ denotes the element of $\cL^2_{-n-2}$ such that 
$$
\widetilde F(z,t,h) =\frac{1}{(2\pi)^{n+1}} \int_{-\infty}^0
\tr\Big(\tau_F(\lambda)P_0\big(e^{\lambda
  h}\sigma_\lambda[z,t]^*-e^\lambda\sigma_\lambda[0,0]^*\big) \Big)\,
|\lambda|^n d\lambda+\widetilde F(0,0,1)\,.
$$
Also, by the definition of reproducing kernel, \eqref{PW-D-eq2}, and
$(z,t,h)=\Psi(\z)$, and writing $K(\z,\cdot)=K_\z$, we have 
$$
\widetilde F(z,t,h) = F(\z)
= \la F, K_{\z} \ra_{\cD_{(m)}} =
\frac{\Gamma(2m-n-1)}{2^{2m-n-1}(2\pi)^{n+1}} 
\int_{-\infty}^0 \tr \big( \tau_F(\lambda) \tau_{K_\z}(\lambda)^*\big)
\, |\lambda|^{2n+1} d\lambda\,+F(\i)\overline{K_\z(\i)}\,.
$$
Since these two equalities hold for all $\tau\in\cL^2_{-n-2}$ we conclude that 
$$
K_\z(\i)=1 \quad \textrm{and } \quad \tau_{K_\z}(\lambda)= C |\lambda|^{-n-1}P_0\big(e^{\lambda
  h}\sigma_\lambda[z,t]-e^\lambda\sigma_\lambda[0,0]\big)\,,
$$
where $C=\textstyle{ \frac{2^{2m-n-1}}{\Gamma(2m-n-1)}}$.
Thus, from Lemma \ref{NEW2}, we obtain
\begin{align*}
& \widetilde K_{(z,t,h)}(w,s,k)\\
& \quad=\frac{C}{(2\pi)^{n+1}}\int_{-\infty}^0
\tr\Big(P_0\big(e^{\lambda
  h}\sigma_\lambda[z,t]-e^{\lambda}\sigma_{\lambda}[0,0]\big)P_0\big(e^{\lambda
  k}\sigma_\lambda[w,s]^*-e^{\lambda}\sigma_{\lambda}
[0,0]^*\big)\Big)\, |\lambda|^{-1}d\lambda+1\,.
\end{align*}
Exploiting \eqref{P0-sigma} the conclusion follows.
\ms
\qed

We conclude the section providing the justification for referring to
the space $\cD_{(m)}$ as the Dirichlet space on the Siegel
half-space. We denote by $\dot{\cD}(B)$ the Dirichlet space
modulo the constant functions on the unit ball
$B\subseteq\bbC^{n+1}$. If $f$ is 
  holomorphic on $B$, $f(\z)=\sum_{|\alpha|\ge0}  a_\alpha
  \z^\alpha$, the norm in $\dot{\cD}(B)$ is given by
$$
\| f\|_{\dot{\cD}(B)}^2 
= \sum_{|\alpha|\ge0} |\alpha|\frac{\alpha!}{|\alpha|!} |a_\alpha|^2\,,
$$
see \cite{Zhu-TAMS}.
Then, the reproducing kernel of $\dot{\cD}(B)$ is given by
$$
K^B(\omega,\z)=
\frac{(n+1)!}{\pi^{n+1}}
\log\frac{1}{1- \omega\cdot\ov\z }\,,
$$
see, e.g. \cite{P} or \cite{Zhu}.

\begin{thm}\label{Mobius-invariance}
 Let $m$ be a positive integer, $m>(n+1)/2$ and denote by
 $\dot{\cD}_{(m)}$ the space $\cD_{(m)}/\bbC$, that is, the space
 $\cD_{(m)}$ modulo the constant functions.  Then: 

{\rm{(1)}} The space $\dot{\cD}_{(m)}$, identified with
$\cD_{(m)}(\i)$  and
endowed with the norm
\begin{equation*}
\|F\|_{\dot{\cD}_{(m)}}= \|\p^m_{\z_{n+1}}F\|_{A^2_{2m-n-2}} \,,
\end{equation*}
is a Hilbert space with reproducing kernel
\begin{equation}\label{kernel-Dir-dot}
K(\omega,\z)=\gamma_{n,m}\log\frac{Q(\omega,\i)Q(\i, \z)}{Q(\omega,\z)}\,,
\end{equation}
where $\gamma_{n,m}$ is as in Corollary \ref{repr-ker--Dir-cor}.

{\rm (2)} For every $\varphi\in \operatorname{Aut}(\cU)$ and every
$F\in\dot{\cD}_{(m)}$ it holds 
\begin{equation}\label{norm-aut}
\|F\circ\varphi\|_{\dot{\cD}_{(m)}}=\|F\|_{\dot{\cD}_{(m)}}\,.
\end{equation}

{\rm (3)} The space $\dot{\cD}_{(m)}$ is isometrically equivalent to
$\dot{\cD}(B)$, the Dirichlet space modulo the constant functions on
the unit ball $B\subset\bbC^{n+1}$. In particular, the space
$\dot{\cD}_{(m)}$ is the unique  Hilbert space of holomorphic
functions on $\cU $ satisfying property \eqref{norm-aut}. 
\end{thm}
\proof
The proof of (1) is straightforward. We now prove (2).
It is enough to prove it for a $\varphi\in\operatorname{Aut}(\cU)$ of
type (i), (ii), (iii) and (iv) described in Lemma \ref{Aut-U}. If
$\varphi$ falls in the cases (i), (ii) or (iii), then it is immediate
to obtain \eqref{norm-aut} by direct computations.  If $\varphi$ is of
type (iv) 
we observe that 
\begin{align*}
(K_{\omega}\circ\varphi)(\z)
&=K(\varphi(\z),\omega)=K(\z,\varphi(\omega))=K_{\varphi(\omega)}(\z)
\,. 
\end{align*}
In particular,
$$
K(\varphi(\z),\varphi(\omega))=K(\z,\omega) \,.
$$
Thus, if $\omega_1,\dots,\omega_N\in\cU$ and
$F(\z)=\sum_{k=1}^N\alpha_{k}K(\z,\omega_k)$, we have
\begin{align*}
\|F\circ\varphi\|^2_{\dot{\cD}_{(m)}}
&=\Big\|\sum_{k=1}^N\alpha_k K(\varphi(\cdot),\omega_k)\Big\|^2_{\dot{\cD}_{(m)}}
=\Big\|\sum_{k=1}^N\alpha_k K(\cdot,\varphi(\omega_k))\Big\|^2_{\dot{\cD}_{(m)}}\\
&=\sum_{j,k=1}^N\overline{\alpha_j}\alpha_k
K(\varphi(\omega_j),\varphi(\omega_k)) 
=\sum_{j,k=1}^N\overline{\alpha_j}\alpha_k K(\omega_j,\omega_k)\\
&=\|F\|^2_{\dot{\cD}_{(m)}} \,,
\end{align*}
as we wished to show. Since the functions of the form
$\sum_{k=1}^N\alpha_{k}K(\z,\omega_k)$ are dense in $\dot{\cD}_{(m)}$ the
conclusion for a generic $F\in\dot{\cD}_{(m)}$ follows. 

At last, (3) follows from the following observation. Let $K^B$
and $K^\cU$ denote the reproducing kernel of $\dot{D}(B)$ and
$\dot{\cD}_{(m)}$ respectively. Then, 
$$
K^B(\omega,\z)= \frac{(n+1)!}{\pi^{n+1}\gamma_{n,m}}K^\cU (\cC(\omega),\cC(\z))\,,
$$
where $\cC$ denote the Cayley transform and $\gamma_{n,m}$ is as in \eqref{kernel-Dir-dot}.
From this it is easy to deduce that the map
$$
F\mapsto \frac{(n+1)!}{\pi^{n+1}\gamma_{n,m}}(F\circ \cC)
$$
is a surjective isometry from $\dot{\cD}_{(m)}$ onto $\dot{\cD}(B)$ as
we wished to show. Hence, the uniqueness of $\dot{\cD}_{(m)}$ follows
from the analogous result for the space $\dot{\cD}(B)$, see
\cite{AF,Zhu-TAMS}. 
\ms
\qed

\section{The Drury--Arveson norm on the unit ball and final remarks}
\label{remarks-sec}

Following \cite{P}, we set
$$
\sR_0 = \text{Id},\qquad \sR_k= \Big( \text{Id} +\frac{R}{k}\Big)
\sR_{k-1}\quad \text{for } k=1,2,\dots,
$$
where $R$ denotes the radial derivative $Rf(\z)=\sum_{j=1}^{n+1} \z_j
\p_{\z_j} f(\z)$. 
Then, we
have the following result on the exact norm in $\DA(B)$. 

\begin{thm}  \label{DA-ball-thm}
  If $f\in\DA(B)$, then
$$
\|f\|_{\DA(B)}^2 
= n\frac{n!}{\pi^{n+1}} \int_B \frac{(1-|\z|^2)^{n-1} }{|\z|^{2n}} \big| \sR_{n}
f(\z)\big|^2 d\z \,.
$$
\end{thm}

This is an elementary computation  that follows from the indenty $\sR_n z^\alpha =
\frac{(n+|\alpha|)!}{n!|\alpha|!} z^\alpha$, for every multiindex
$\alpha$. 
\ms

We believe this work raises some interesting questions.  We first mention 
the characterization of the Carleson
measures and of the multiplier algebra for the scale of spaces studied
in this work.  Moreover, these spaces  depend on the parameter 
$\nu$ where $\nu\geq -n-2$.  It would be interesting to study the
spaces corresponding to the values  $\nu<-n-2$. 
Furthermore, we would like to study the analogous Banach spaces, whose
underlying norm is the $L^p$-norm, with  $p\neq2$. 
Finally, it would be interesting to extend the results in this paper
to the more general setting of Siegel domains of
type $I\!I$ . We plan to come back to these problems in
future works.

\bibliography{DA-bib}

\newcommand{\etalchar}[1]{$^{#1}$}
\providecommand{\bysame}{\leavevmode\hbox to3em{\hrulefill}\thinspace}
\providecommand{\MR}{\relax\ifhmode\unskip\space\fi MR }
\providecommand{\MRhref}[2]{%
  \href{http://www.ams.org/mathscinet-getitem?mr=#1}{#2}
}
\providecommand{\href}[2]{#2}
\begin{thebibliography}{DGGMR07}

\bibitem[AF85]{AF}
J.~Arazy and S.~D. Fisher, \emph{The uniqueness of the {D}irichlet space among
  {M}\"obius-invariant {H}ilbert spaces}, Illinois J. Math. \textbf{29} (1985),
  no.~3, 449--462. \MR{786732}

\bibitem[ARS08]{ARS}
N.~Arcozzi, R.~Rochberg, and E.~Sawyer, \emph{Carleson measures for the
  {D}rury-{A}rveson {H}ardy space and other {B}esov-{S}obolev spaces on complex
  balls}, Adv. Math. \textbf{218} (2008), no.~4, 1107--1180. \MR{2419381}

\bibitem[ARS10]{ARSW-2variations}
\bysame, \emph{Two variations on the {D}rury-{A}rveson space}, Hilbert spaces
  of analytic functions, CRM Proc. Lecture Notes, vol.~51, Amer. Math. Soc.,
  Providence, RI, 2010, pp.~41--58. \MR{2648865}

\bibitem[Arv72]{Arveson}
W.~Arveson, \emph{Subalgebras of {$C^{\ast} $}-algebras. {II}}, Acta Math.
  \textbf{128} (1972), no.~3-4, 271--308. \MR{0394232}

\bibitem[BB89]{BB}
F.~Beatrous and J.~Burbea, \emph{Holomorphic {S}obolev spaces on the ball},
  Dissertationes Math. (Rozprawy Mat.) \textbf{276} (1989), 60. \MR{1010151}

\bibitem[BBG{\etalchar{+}}04]{camerun}
D.~B\'ekoll\'e, A.~Bonami, G.~Garrig\'os, C.~Nana, M.~M. Peloso, and F.~Ricci,
  \emph{Lecture notes on {B}ergman projectors in tube domains over cones: an
  analytic and geometric viewpoint}, IMHOTEP J. Afr. Math. Pures Appl.
  \textbf{5} (2004), Exp. I, front matter + ii + 75. \MR{2244169}

\bibitem[CSW11]{CSW}
{\c S}.~Costea, E.~Sawyer, and B.~D. Wick, \emph{The corona theorem for the
  {D}rury-{A}rveson {H}ardy space and other holomorphic {B}esov-{S}obolev
  spaces on the unit ball in {$\Bbb C^n$}}, Anal. PDE \textbf{4} (2011), no.~4,
  499--550. \MR{3077143}

\bibitem[DGGMR07]{DG}
P.~Duren, E.~A. Gallardo-Guti\'errez, and A.~Montes-Rodriguez, \emph{A
  {P}aley-{W}iener theorem for {B}ergman spaces with application to invariant
  subspaces}, Bull. Lond. Math. Soc. \textbf{39} (2007), no.~3, 459--466.
  \MR{2331575}

\bibitem[Dru78]{Drury}
S.~W. Drury, \emph{A generalization of von {N}eumann's inequality to the
  complex ball}, Proc. Amer. Math. Soc. \textbf{68} (1978), no.~3, 300--304.
  \MR{480362}

\bibitem[Fel91]{Feldman}
M.~Feldman, \emph{Mean oscillation, weighted {B}ergman spaces, and {B}esov
  spaces on the {H}eisenberg group and atomic decompositions}, J. Math. Anal.
  Appl. \textbf{158} (1991), no.~2, 376--395. \MR{1117569}

\bibitem[FK94]{FK}
J.~Faraut and A.~Kor\'anyi, \emph{Analysis on symmetric cones}, Oxford
  Mathematical Monographs, The Clarendon Press, Oxford University Press, New
  York, 1994, Oxford Science Publications. \MR{1446489}

\bibitem[Fol89]{Folland}
G.~B. Folland, \emph{Harmonic analysis in phase space}, Annals of Mathematics
  Studies, vol. 122, Princeton University Press, Princeton, NJ, 1989.
  \MR{983366}

\bibitem[Kuc17a]{Kucik}
A.~S. Kucik, \emph{Carleson measures for {H}ilbert spaces of analytic functions
  on the complex half-plane}, J. Math. Anal. Appl. \textbf{445} (2017), no.~1,
  476--497. \MR{3543778}

\bibitem[Kuc17b]{Kucik2}
\bysame, \emph{Multipliers of {H}ilbert spaces of analytic functions on the
  complex half-plane}, Oper. Matrices \textbf{11} (2017), no.~2, 435--453.
  \MR{3655659}

\bibitem[Lax59]{Lax}
P.~D. Lax, \emph{Translation invariant spaces}, Acta Math. \textbf{101} (1959),
  163--178. \MR{0105620}

\bibitem[OV79]{OV}
R.~D. Ogden and S.~V\'agi, \emph{Harmonic analysis of a nilpotent group and
  function theory of {S}iegel domains of type {${\rm II}$}}, Adv. in Math.
  \textbf{33} (1979), no.~1, 31--92. \MR{540636}

\bibitem[Pel92]{P}
M.~M. Peloso, \emph{M\"obius invariant spaces on the unit ball}, Michigan Math.
  J. \textbf{39} (1992), no.~3, 509--536. \MR{1182505}

\bibitem[Pol15]{Poltoratski}
A.~Poltoratski, \emph{Toeplitz approach to problems of the uncertainty
  principle}, CBMS Regional Conference Series in Mathematics, vol. 121,
  Published for the Conference Board of the Mathematical Sciences, Washington,
  DC; by the American Mathematical Society, Providence, RI, 2015. \MR{3309830}

\bibitem[PW34]{PW}
R.~E. A.~C. Paley and N.~Wiener, \emph{Fourier transforms in the complex
  domain}, American Mathematical Society Colloquium Publications, vol.~19,
  American Mathematical Society, Providence, RI, 1987[1934], Reprint of the
  1934 original. \MR{1451142}

\bibitem[Ric92]{R}
F.~Ricci, \emph{Harmonic {A}nalysis on the {H}eisenberg {G}roup}, Unpublished
  notes, Cortona (1992), ii+67.

\bibitem[RS16]{RS}
S.~Richter and J.~Sunkes, \emph{Hankel operators, invariant subspaces, and
  cyclic vectors in the {D}rury-{A}rveson space}, Proc. Amer. Math. Soc.
  \textbf{144} (2016), no.~6, 2575--2586. \MR{3477074}

\bibitem[RT83]{RT}
F.~Ricci and M.~Taibleson, \emph{Boundary values of harmonic functions in mixed
  norm spaces and their atomic structure}, Ann. Scuola Norm. Sup. Pisa Cl. Sci.
  (4) \textbf{10} (1983), no.~1, 1--54. \MR{713108}

\bibitem[Sei04]{Seip}
K.~Seip, \emph{Interpolation and sampling in spaces of analytic functions},
  University Lecture Series, vol.~33, American Mathematical Society,
  Providence, RI, 2004. \MR{2040080}

\bibitem[Ste93]{Stein}
E.~M. Stein, \emph{Harmonic analysis: real-variable methods, orthogonality, and
  oscillatory integrals}, Princeton Mathematical Series, vol.~43, Princeton
  University Press, Princeton, NJ, 1993, With the assistance of Timothy S.
  Murphy, Monographs in Harmonic Analysis, III. \MR{1232192}

\bibitem[Tch08]{T}
E.~Tchoundja, \emph{Carleson measures for the generalized {B}ergman spaces via
  a {$T(1)$}-type theorem}, Ark. Mat. \textbf{46} (2008), no.~2, 377--406.
  \MR{2430733}

\bibitem[VW12]{VW}
A.~Volberg and B.~D. Wick, \emph{Bergman-type singular integral operators and
  the characterization of {C}arleson measures for {B}esov-{S}obolev spaces and
  the complex ball}, Amer. J. Math. \textbf{134} (2012), no.~4, 949--992.
  \MR{2956255}

\bibitem[Wal73]{Wallach}
N.~R. Wallach, \emph{Harmonic analysis on homogeneous spaces}, Marcel Dekker,
  Inc., New York, 1973, Pure and Applied Mathematics, No. 19. \MR{0498996}

\bibitem[Zhu91]{Zhu-TAMS}
K.~Zhu, \emph{M\"obius invariant {H}ilbert spaces of holomorphic functions in
  the unit ball of {${\bf C}^n$}}, Trans. Amer. Math. Soc. \textbf{323} (1991),
  no.~2, 823--842. \MR{982233}

\bibitem[Zhu05]{Zhu}
\bysame, \emph{Spaces of holomorphic functions in the unit ball}, Graduate
  Texts in Mathematics, vol. 226, Springer-Verlag, New York, 2005. \MR{2115155}

\end{thebibliography}
\bibliographystyle{amsalpha}

\end{document}